\documentclass[11pt, a4paper]{amsart}
\usepackage{latexsym}
\usepackage{graphicx,amsmath,amssymb,cases,hyperref,float,enumitem,mathrsfs} 
\usepackage{color}
\usepackage{tikz}
\usepackage{tikz-cd}
\usepackage{array,booktabs,longtable}
\usepackage{verbatim}
\usepackage[dvipsnames]{xcolor}
\usepackage{ytableau}
\usepackage[all]{xy}
\SelectTips{cm}{12}
\objectmargin={1pt}
\newsavebox{\AFGDTableBox}
\newsavebox{\AFGQTableBox}
\usepackage[
  a4paper,
  left=30mm,
  right=30mm,
  top=25mm,
  bottom=25mm
]{geometry} 
\usepackage{geometry}

\usepackage{graphicx}
\usepackage{booktabs}
\usepackage{pdflscape}

\usepackage[mathscr]{euscript}

\newcommand{\ZZ}{\mathbb Z}

\newcommand{\Waf}{W_{\mathrm{af}}}
\newcommand{\Wgr}{W_{\mathrm{af}}^{0}}
\newcommand{\gh}{ k}
\newcommand{\cg}{\ell}
\newcommand{\ck}{\ell}
\newcommand{\kh}{k}
\newcommand{\eps}{\varepsilon}
\newcommand{\id}{\mathrm{id}}

\newcommand{\Supp}{\operatorname{Supp}}

\newcommand{\Gr}{\mathrm{Gr}}
\newcommand{\af}{\mathrm{af}}

 \newcommand{\CC}{\widetilde C}
\newcommand{\Z}{\mathbb Z}
\newcommand{\W}{W_{\rm af}}
\newcommand{\Wg}{W_{\rm af}^{0}}

\title[Equivariant $K$-homology Pieri rules]{Pieri rules in the torus-equivariant $K$-homology \\ of affine Grassmannians of type $A$}
\author[K. Brahma]{Koushik Brahma}
\address{
Waseda University, Faculty of Science and Engineering,
3-4-1 Okubo, Shinjuku-ku,
Tokyo 169-8555, Japan}
\email{koushikbrahma95@gmail.com}
\date{\today}

\newtheorem{innercustomgeneric}{\customgenericname}
\providecommand{\customgenericname}{}
\newcommand{\newcustomtheorem}[2]{%
	\newenvironment{#1}[1]
	{%
		\renewcommand\customgenericname{#2}%
		\renewcommand\theinnercustomgeneric{##1}%
		\innercustomgeneric
	}
	{\endinnercustomgeneric}
}

\newcustomtheorem{customthm}{Theorem}
\newcustomtheorem{customcor}{Corollary}
\newcustomtheorem{customprop}{Proposition}

\newtheorem{theorem}{Theorem}[section]
\newtheorem{lemma}[theorem]{Lemma}
\newtheorem{proposition}[theorem]{Proposition}
\newtheorem{corollary}[theorem]{Corollary}
\newtheorem{conjecture}[theorem]{Conjecture}
\theoremstyle{remark}
\newtheorem{definition}[theorem]{Definition}
\newtheorem{example}[theorem]{Example}

\newtheorem{remark}[theorem]{Remark}
\theoremstyle{plain}

\numberwithin{equation}{section}

\theoremstyle{remark}

\begin{document}


\begin{abstract}

In this paper, we give explicit Pieri rules in the Schubert structure-sheaf and boundary ideal-sheaf bases of the torus-equivariant $K$-homology of the affine Grassmannians of $SL_{n+1}$.
We express both Pieri classes explicitly in terms of translation elements on signed coordinate characters using the star action on the $K$-theoretic Peterson algebra and polynomial interpolation. A classification of the elements below these translations in Bruhat order yields finite formulas for their coefficients in the affine $K$-nilHecke algebra. We prove signed positivity for the
equivariant Pieri structure constants in the boundary ideal-sheaf basis. In non-equivariant $K$-homology, we obtain explicit binomial formulas for the structure constants in both Schubert bases.
\end{abstract}

\maketitle

\section{Introduction}
Let $G$ be a simple simply connected complex algebraic group and $T\subset G$ be the maximal torus.
The affine Grassmannian corresponding to $G$ is the ind-projective variety
\[
 \mathrm{Gr}_G=G(\mathbb C((z)))/G(\mathbb C[[z]]).
\]
The $T$-equivariant $K$-homology $K^T_*(\Gr_G)$ of $\Gr_G$ is the Grothendieck group of $T$-equivariant coherent sheaves on $\Gr_G$ with proper support. 
The based-loop-space model of $\Gr_G$ induces a Pontryagin product on $K^T_*(\Gr_G)$, making $K^T_*(\Gr_G)$ into a commutative $R(T)$-algebra. The readers are referred to \cite{LSS,ISY} for more details on $T$-equivariant $K$-homology of affine Grassmannian.

Let $W$ denote the finite Weyl group of $G$. Its affine Weyl group $\W$ is the semidirect product $\W=Q^\vee\rtimes W$, where $Q^\vee$ is the coroot lattice. The affine Grassmannian elements $\Wgr$ are the minimal-length representatives of the cosets in $\W/W$. For $w\in W_{\mathrm{af}}^0$, the Schubert cell indexed by $w$ is the Iwahori orbit of the corresponding $T$-fixed point in $\mathrm{Gr}_G$. Let $X_w$ be the closure of the Schubert cell indexed by $w$, and $\partial X_w$ for its reduced boundary. Then $X_u\subseteq X_w$ if and only if $u\le w$ in the Bruhat order. The Schubert structure-sheaf basis and boundary ideal-sheaf basis are given by
\[
\mathcal O_w=[\mathcal O_{X_w}], \qquad \xi_w=[\mathcal I_{\partial X_w}].
\]
Here $\mathcal{O}_{X_w}$ is the structure sheaf of $X_w$, and $\mathcal{I}_{\partial X_w}$ is the ideal sheaf of the boundary $\partial X_w$ in $X_w$.

 The study of Schubert calculus in the \(T\)-equivariant homology \(H_*^T(\mathrm{Gr}_G)\) of the affine Grassmannian is motivated by Peterson’s “quantum equals affine” philosophy \cite{Pet} (see also \cite{LS-Q}). After suitable localization, the Peterson isomorphism identifies \(H_*^T(\mathrm{Gr}_G)\) with the \(T\)-equivariant quantum cohomology \(QH_T^*(G/B)\) of the flag variety. The correspondence between Schubert classes under this isomorphism expresses structure constants in \(QH_T^*(G/B)\) in terms of structure constants in \(H_*^T(\mathrm{Gr}_G)\). The \(K\)-theoretic analogue of Peterson isomorphism was studied in the non-equivariant setting in type \(A\) in \cite{IIM}. The correspondence in general type was conjectured in \cite{LLMSC} and proved by Kato \cite{Kato}. Chow and Leung \cite{CL} also gave an independent proof.

Pieri rules for the homology of affine Grassmannians were established by Lam, Lapointe, Morse and Shimozono \cite{LLMS} in type~A, by Lam, Schilling and Shimozono \cite{LSSC} in type~C, and by Pon \cite{Pon} in types~B and D. In $T$-equivariant homology, Lam and Shimozono \cite{LS} established a Pieri rule in type~A and conjectured analogous rules in types B and C. In non-equivariant \(K\)-homology of type \(A\), Lam, Schilling and Shimozono \cite{LSS} describe the ideal-sheaf Pieri elements in the $K$-nilHecke algebra and deduce the Pieri rule in ideal-sheaf basis. Takigiku \cite{Tak} gives a Pieri rule and a rectangle factorization formula for sums of \(K\)-theoretic \(k\)-Schur functions, which represent the non-equivariant structure-sheaf.
Ikeda, Shimozono, and Yamaguchi \cite{ISY} construct a symmetric-function realization of both Schubert bases in $T$-equivariant $K$-homology of affine Grassmannian in type A.

For $G=SL_{n+1}$, we consider the Pieri elements \(\rho_r\) for \(1\le r\le n\), and compute the expansion of the Peterson representative \(k_{\rho_r}\) of the boundary ideal-sheaf class \(\xi_{\rho_r}\) in the \(T_w\)-basis as well as the expansion of the Peterson representative \(\ell_{\rho_r}\) of the structure-sheaf class \(\mathcal O_{\rho_r}\) in the \(D_w\)-basis. These expansions in the level-zero affine \(K\)-nilHecke algebra give finite-sum formulas for the  Pieri structure constants in both Schubert bases. To compute the coefficients, we attach a combinatorial data on the support of the affine Weyl group elements and use the nilHecke Yang--Baxter operators of Fomin and Kirillov \cite{FK1, FK2}, see also Lascoux, Leclerc and Thibon \cite{LLT}. We also use the subword criterion, the lifting property, and parabolic factorization for
Coxeter groups~\cite{BB}.

Signed positivity for Schubert structure constants in equivariant
$K$-theory was established by Anderson, Griffeth, and Miller
\cite{AGM2011} for finite-dimensional flag varieties and by Baldwin
and Kumar \cite{BK2017} for opposite Schubert structure sheaves on
thick Kac--Moody flag varieties. For the Pontryagin product in the
equivariant $K$-homology of the affine Grassmannian, Kumar
\cite{Kumar} conjectured signed positivity of the structure
constants in the Schubert structure-sheaf basis. We establish signed positivity for the Pieri structure constants in the boundary ideal-sheaf basis for the $T$-equivariant $K$-homology of affine Grassmannians of $SL_{n+1}$.

\subsection{The level-zero affine \texorpdfstring{$K$}{K}-nilHecke algebra}
Kostant and Kumar introduced the nilHecke ring in equivariant cohomology \cite{KK} and developed its \(K\)-theoretic analogue in \cite{KK1}. We use the level-zero affine analogue following \cite[Section 2]{LSS} and \cite[Section 2]{ISY}. 

The finite Weyl group $W$ acts naturally on $R(T)$. 
The level-zero action of the affine Weyl group $\W$ on $R(T)$ is obtained via the group homomorphism $\W\to W$, given by $t_\lambda u\mapsto u$ for $\lambda\in Q^\lor$ and $u\in W$. Let $\Phi$ be the finite root system and let $\Delta$ be the multiplicative subset of $R(T)$ generated by $\{1-e^\alpha\mid \alpha\in\Phi\}$. The level-zero action extends to the localization $R(T)^\Delta=\Delta^{-1}R(T)$. We write $\mathbb K^\Delta=R(T)^\Delta\otimes_\Z\Z[\W]$ for the twisted group algebra with multiplication
\begin{equation}\label{eq_tw_act}
  (fw)(gv)=f\,w(g)\,wv,
\qquad f,g\in R(T)^\Delta,\quad v,w\in\W.  
\end{equation}

Let $\{\alpha_i~|~1\le i\le n\}$ be the simple roots in $\Phi$ and $\alpha_0=-\theta$, where $\theta$ is the highest root. For $i\in \{0,1,\dots,n\}$, define
\begin{equation}\label{eq:convention}
 T_i=(1-e^{\alpha_i})^{-1}(s_i-1), \quad D_i=T_i+1.
\end{equation} 
Direct multiplication gives
\begin{equation}\label{eq:cross}
 T_i^2=-T_i,\qquad D_i^2=D_i.
\end{equation}
The elements $T_i$ and $D_i$ satisfy the Coxeter braid relations.
If $w \in W_\af$, one can define
\begin{equation}\label{eq_pro}
    T_w:= T_{i_1} \cdots T_{i_\ell},~ D_w:= D_{i_1} \cdots D_{i_\ell} \text{ for a reduced expression $w = s_{i_1} \cdots s_{i_\ell}$.}
\end{equation} 
The level-zero affine $K$-nilHecke algebra $\mathbb K$ is defined to be the left $R(T)$-module generated by $\{T_w:w\in \W\}$ as well as $\{D_w:w\in \W\}$. More precisely,
$$\mathbb K=\bigoplus_{w\in\W}R(T)T_w=\bigoplus_{w\in\W}R(T)D_w.$$
Then $\mathbb K$ has a commutative subalgebra $\mathbb L$ known as the $K$-theoretic Peterson algebra which is defined as the centralizer of $R(T)$ in $\mathbb K$. 
 Moreover, there is a Hopf algebra isomorphism 
 $$K^T_*( \text{Gr}_G ) \cong \mathbb L,$$ see~\cite[Theorem 5.3]{LSS} as an $R(T)$-algebra. 
Let $k_w$ and $\ell_w$ denote the images of the boundary ideal-sheaf $\xi_w$ and structure-sheaf $\mathcal{O}_w$, respectively. 
Then we have the following
$$\mathbb L \cong \bigoplus_{w\in \Wgr} R(T)\,k_w=\bigoplus_{w\in \Wgr} R(T)\,\ell_w.$$

Let $\xi_w^0$ and \(\mathcal{O}_{w}^0\) denote the corresponding non-equivariant Schubert classes in $K_*(\Gr_{G})$.
Let $\phi_0\colon R(T)\to\mathbb Z$ be the augmentation homomorphism
defined by $\phi_0(e^\lambda)=1$ for every $\lambda\in X^*(T)$ (the character lattice).
We recall the $0$-Hecke ring \cite[Section 2.3]{LSS} $$\mathbb K_0=\bigoplus_{w\in \W}\Z\, T_w.$$ Then $\phi_0$ induces a map $$\phi_0:\mathbb K\to \mathbb K_0, \quad \phi_0\big(\sum_{w\in \W}a_w T_w\big)=\sum_{w\in \W}\phi_0(a_w)T_w.$$ Then $\mathbb K_0$ has a subalgebra $\mathbb L_0$ (known as 'the $K$-affine Fomin-Stanley algebra') defined by $$\mathbb L_0=\{b\in \mathbb K_0~|~\phi_0(bq)=\phi_0(q)b \text{ for all } q\in R(T)\}.$$
Moreover, there is a Hopf-algebra isomorphism 
$$K_*(\Gr_G)\cong \mathbb L_0,$$
sending $\xi_w^0$ and $\mathcal O_w^0$ to $\phi_0(\gh_w)$ and $\phi_0(\ck_w)$ respectively, see \cite[Theorem 6.8]{LSS}.

\subsection{Pieri classes and Pieri structure constants}
In type $\mathrm{A}_n$, we consider the Pieri elements (single-row affine Grassmannian elements) 
\[
 \rho_r=s_{r-1}\cdots s_1s_0\quad(1\le r\le n),
 \qquad \rho_0=\id,
\]
see \cite[Section 1.2]{LLMSS}. Set $\gh_r=\kh_{\rho_r}$ and $\ck_r=\ck_{\rho_r}$ for $1\le r\le n$ with $\gh_0=\ck_0=1$, call $\kh_r$ and $\ck_r$ as the ordinary Pieri classes and closed Pieri classes respectively. Let $I=\{0,1,\ldots,n\}$. In the level-zero affine $K$-nilHecke algebra, write 
\begin{equation}\label{eq:intromain}
 \gh_r=\sum_{w\in W_{\rm af}}C_{r,w}T_w, \qquad \ck_r=\sum_{w\in W_{\rm af}}E_{r,w}D_w.
\end{equation}

We give finite formulas for $C_{r,w}$ and $E_{r,w}$ in terms of the
`admissible support data' $(S,L,R)$ of $w$, where $S=\Supp(w)$ and
$L,R\subseteq S$, see Section~\ref{sec:colored}. Let $c$ be the number
of connected components of $S$ for proper support, and $c=0$ for
full support. 

Fix $1\le r\le n$. Define $Z:=R\setminus\{0\}$, $A:=Z\cap[1,r]$, and
\[
 [1,r]\setminus A=\{j_1<\cdots<j_k\},\qquad
 h=c-1+\mathbf{1}_{0\in R},\qquad a=n-c-|R|,
\]
where $\mathbf{1}_{0\in R}$ is $1$ if $0\in R$ and $0$ otherwise.
Writing $Z\setminus A=\{i_1<\cdots<i_s\}$, set
\begin{equation}\label{eq:jp-list}
 (b_1,\ldots,b_d)
 =(\underbrace{0,\ldots,0}_{h\text{ copies}},i_1,\ldots,i_s).
\end{equation}

For a set $M$ of positive integers, write
$M_\ell=1+|M\cap[1,\ell-1]|$.
For $1\le i\le n$, set $x_i=e^{-\alpha_i}$ and
$y_i=x_i\cdots x_n$. We use the conventions $y_0=y_{n+1}=1$ and set
$Y_r=\prod_{i=1}^r y_i$ for $0\le r\le n$, with $Y_0=1$.
Consider the polynomial semiring
\[
 \mathcal P=\mathbb Z_{\ge0}[x_1-1,\ldots,x_n-1].
\]
\begin{customthm}{A}[Theorem \ref{thm:jp-ordinary}, Theorem \ref{thm:jp-closed}]
Let $w\in\W\setminus\{\mathrm{id}\}$ and fix $1\le r\le n$. If $w$
has no admissible support datum or, its datum gives $k=0$ then $C_{r,w}=E_{r,w}=0$. If $k\ge1$, the following
formulas hold.
\[
\begin{aligned}
(-1)^{\ell(w)-r}C_{r,w}
&=\frac{Y_{r-1}\prod_{i\in L}(y_i-1)}{\prod_{i\in A}y_i} \sum_{\substack{M\subseteq[1,d]\\|M|=k-1}}
 \left(\prod_{\ell\in M}\frac{1}
 {y_{b_\ell}}\right)
 \prod_{\ell\notin M}
 \left(\frac{y_{j_{M_\ell}}}{y_{b_\ell}}-1\right).\\
(-1)^{\ell(w)-r}E_{r,w}
&=\frac{Y_r\prod_{i\in L}(y_i-1)}
{\prod_{i\in L}y_i}\sum_{\substack{M\subseteq[1,a+d]\\|M|=k-1}}
\left(\prod_{\substack{1\le\ell\le a\\\ell\notin M}}
y_{j_{M_\ell}}\right)
\prod_{\substack{1\le t\le d\\a+t\notin M}}
\left(y_{j_{M_{a+t}}}-y_{b_t}\right).
\end{aligned}
\] 
In particular, $C_{r,w}=E_{r,w}=0$ whenever $\ell(w)>2n$. For the identity element,
\[
 C_{r,\mathrm{id}}=0,\qquad
 (-1)^rE_{r,\mathrm{id}}
   =Y_r\sum_{s=0}^{n-r}\mathsf h_s(y_1,\ldots,y_r),
\]
where $\mathsf h_s$ is the complete homogeneous symmetric polynomial of degree
$s$. 
\end{customthm}

\begin{customcor}{B}
   For all $w\in \W$, $(-1)^{\ell(w)-r}C_{r,w}$ and $(-1)^{\ell(w)-r}E_{r,w}$ both lies in $\mathcal P.$
\end{customcor}

For $v\in\Wg$ and $1\le r\le n$, consider the multiplication
\begin{equation}\label{eq_st_con}
   \gh_r\kh_v=\sum_{u\in\Wg}d_{r,v}^{\,u}\kh_u,\qquad \cg_r\ck_v=\sum_{u\in\Wg}\widetilde d_{r,v}^{\,u}\ck_u. 
\end{equation}

The structure constants $d_{r,v}^{\,u}$ and $\widetilde d_{r,v}^{\,u}$ can be written explicitly as a finite sum of $C_{r,w}$ and $E_{r,w}$ respectively for certain $w\in \W$, see Proposition~\ref{st_con}.

\begin{customthm}{C}[Theorem \ref{thm:structure}]
For all $u,v\in\Wgr$ and $1\le r\le n$, the sign-adjusted Pieri coefficient
$(-1)^{\ell(u)-\ell(v)-r}d_{r,v}^{\,u}$ belongs to $\mathcal P$.
\end{customthm}

For a proper subset \(S\subsetneq I\), define \(d_S\in W_{\mathrm{af}}\) to be the element represented by a word containing \(s_i\) for each \(i\in S\), exactly once, with \(s_i\) occurring to the left of \(s_{i-1}\), whenever \(i,i-1\in S\). Indices are interpreted modulo \(n+1\). Set \(d_\varnothing=\mathrm{id}\).
\begin{customthm}{D}
For every $1\leq r\leq n$, we have
\begin{align}
     \phi_0(\gh_r)&=\sum_{\substack{S\subsetneq I\\|S|=r}}T_{d_S}.\nonumber\\
\phi_0(\ck_r)&=
\sum_{\substack{S\subsetneq I\\|S|\le r}}
(-1)^{r-|S|}
\binom{n-|S|}{r-|S|}
D_{d_S}\nonumber.
\end{align}    
\end{customthm}
The cyclically decreasing expansion of \(\phi_0(k_r)\) is due to Lam, Schilling and Shimozono \cite[Corollary 7.5]{LSS}. However, we derive the expansions of \(\phi_0(k_r)\) and \(\phi_0(\ell_r)\) directly by applying $\phi_0$ to our equivariant coefficient formulas.
 As a corollary, we have
    $$(-1)^{\ell(u)-\ell(v)-r}\phi_0(d_{r,v}^{\,u})\in \Z_{\geq 0} \text{ and }(-1)^{\ell(u)-\ell(v)-r}\phi_0(\widetilde{d}_{r,v}^{\,u})\in \Z_{\geq 0},$$ 
for every $u,v\in \Wgr$ and $1\leq r\le n$. Moreover, we derive an explicit single signed binomial formula for $\phi_0(d_{r,v}^{\,u})$ and $\phi_0(\widetilde{d}_{r,v}^{\,u})$ in Theorem \ref{thm:ns-ideal}.

\subsection{Organization}
We fix the algebraic conventions and identifies the ordinary and closed Pieri classes directly from the star action and interpolation in \S\ref{sec:algebra}. In \S\ref{sec_aux_alg}, we explore the auxiliary algebra which simplifies the computation of the coefficients of the ordinary and closed Pieri classes.
In \S\ref{sec:colored}, we construct a combinatorial data to classify all the affine Weyl group elements those arise in the expansions of $\gh_r$ and $\ck_r$ for $1\le r\le n$. In \S\ref{sec:marks}, we derive the explicit formula of $C_{r,w}$ and $E_{r,w}$ and prove their positivity.
In \S\ref{sec:forget}, we describe both the Pieri classes in non-equivariant $K$-homology. In \S\ref{sec:pieri}, we discuss computing the Pieri structure constants explicitly and prove the positivity of the Pieri structure constants in ideal-sheaf basis. In non-equivariant setting we derive explicit binomial formulas and the positivity for both the structure constants. Some explicit examples are given in \S\ref{sec_ex}.

\section{ The star action on Peterson algebra and interpolation formulas
for Pieri classes}\label{sec:algebra}
In this section, we fix the conventions on roots and characters. We recall star action and identify all the Pieri classes in terms of translation elements. 

\subsection{Conventions on roots and characters}
Let $N:=n+1$ and $I=\{0,1,\ldots,n\}$ denotes the nodes of affine Dynkin diagram. For $n\ge2$, the affine Dynkin diagram is 
\begin{center}
\begin{tikzpicture}[scale=.9, every node/.style={font=\small}]
 \node (a) at (0,0) [circle,draw,inner sep=2pt,label=below:$1$] {};
 \node (b) at (1.3,0) [circle,draw,inner sep=2pt,label=below:$2$] {};
 \node (dots) at (2.6,0) {$\cdots$};
 \node (c) at (3.9,0) [circle,draw,inner sep=2pt,label=below:$n-1$] {};
 \node (d) at (5.2,0) [circle,draw,inner sep=2pt,label=below:$n$] {};
 \node (zero) at (2.6,1.25) [circle,draw,inner sep=2pt,label=above:$0$] {};
 \draw (a)--(b)--(dots)--(c)--(d);\draw (a)--(zero)--(d);
\end{tikzpicture}
\end{center}
We read the indices cyclically by considering modulo $N$. Let $T$ be the diagonal torus of $SL_N$, and let $a_i$
denote its $i$th coordinate character. The simple roots and coroots are
\[
\alpha_i=a_i-a_{i+1},\quad \alpha_i^\vee=\eps_i-\eps_{i+1}\quad (1\le i\le n).
\]
The simple roots form a $\mathbb{Z}$-basis of the root lattice $Q$, and the simple coroots form a $\mathbb{Z}$-basis of the coroot lattice $Q^\vee$. The root system is given by
\[
\Phi=\{\pm(a_i-a_j):1\le i<j\le N\}, \quad \Phi^\vee=\{\pm(\eps_i-\eps_j):1\le i<j\le N\}.
\]

The positive roots are
$a_i-a_j=\alpha_i+\cdots+\alpha_{j-1}$, $i<j$. The highest root is $\theta=a_1-a_N$, with the associated coroot $\theta^\vee=\eps_1-\eps_N$. At level zero, the affine
simple root $\alpha_0=\delta-\theta$ restricts to $\alpha_0=-\theta$, where $\delta$ is the null root.

For $n=1$ the affine Weyl group is the infinite dihedral group. For $n\ge 2$
the affine Weyl group $\W$ is the Coxeter group
\begin{equation*}
    \W=\langle s_0,s_1,\dots,s_n:s_i^2=1,\ (s_is_{i+1})^3=1,
           \ s_is_j=s_js_i\text{ for nonadjacent }i,j\rangle,
\end{equation*}
where $s_{n+1}=s_0$ as we consider the indices modulo $N$. 
The finite Weyl group $W$ is the subgroup of $\W$ given by $$W=\langle s_1,\dots,s_n\rangle\cong S_N.$$

There is an injective group homomorphism $Q^\lor \to \W$ denoted $\mu\mapsto t_\mu$ and the conjugation action of $W$ on $Q^\lor$ is given by $w t_\mu w^{-1} = t_{w(\mu)}$ for $w\in W$ and $\mu\in Q^\lor$.  The affine Weyl group $\W$ is also given by $\W=Q^\vee\rtimes W$. The affine simple
reflection is $$s_0=t_{\theta^\vee}s_\theta=s_{\theta}t_{-\theta^\vee}.$$

If $\ell:\W\to \ZZ_{\ge 0}$ denotes the Coxeter length, the affine Grassmannian elements are
\[
\Wgr=\{w\in \Waf\mid \ell(ws_i)>\ell(w)
\text{ for all $1\le i\le n$}\}.
\]

The Pieri elements are $\rho_r=s_{r-1}\cdots s_1s_0 \in \Wgr$ for all $1\le r\le n$ with $\rho_0=\mathrm{id}$.

In the representation ring $R(T)=\mathbb Z[e^{\pm a_1},\ldots,e^{\pm a_N}]
       /(e^{a_1+\dots+a_N}-1)$,
we set
\[
 x_i=e^{-\alpha_i} \text{ and }
 y_i=e^{-(a_i-a_N)}=x_i\cdots x_n \text{ for all } 1\leq i\le n.
\]
For any positive root $\beta=a_i-a_j$ for $i<j$, we have $e^{-\beta}=y_i/y_j$.
We set $y_0=y_N=1$.

\subsection{The star action of \texorpdfstring{$\mathbb K$}{\mathbb K} on \texorpdfstring{$\mathbb L$}{\mathbb L}}
We recall the star action 
 $$\ast:\mathbb K\times\mathbb L\to\mathbb L$$ following \cite[Section 2.3]{ISY}. Consider $t_\lambda u\in \W$ for $\lambda\in Q^\vee, u\in W$ and $b\in \mathbb L$. Then
 \begin{equation}
(t_\lambda u)\ast b=t_\lambda u b u^{-1}.
\label{eq:star}
\end{equation}
  In particular, for every $i\in I$ 
 $$s_i\ast b= s_i b s_i \text{ for } 1\leq i\le n,\quad  s_0\ast b=t_{\eps_1-\eps_N} s_\theta b s_\theta.$$
For every $i\in I$ and $b\in \mathbb L$, \eqref{eq:convention} gives
$$T_i\ast b=(1-e^{\alpha_i})^{-1}(s_i\ast b-b).$$

The star action induces an $R(T)$-module homomorphism 
$\kappa: \mathbb K\to \mathbb L$ defined by
$\kappa(b)=b\ast 1$ for $b\in \mathbb K.$
The restriction of $\kappa$ to $\mathbb L$ is the identity. Moreover, 
\begin{equation}\label{eq:moduleprojection}
 \kappa(ab)=a\kappa(b)\quad(b\in\mathbb K,~a\in\mathbb L).
\end{equation}
Every element $w\in \W$ has
a unique length-additive factorization $w=w^0v$ with $v\in W$ and $w^0\in\Wg$, see \cite[Proposition~2.4.4]{BB}. For every $w\in \W$, the map $\kappa$ satisfies
\begin{equation}\label{eq:kappa}
\kappa(T_w)=\begin{cases}\kh_w&w\in\Wg\\0&w\notin\Wg\end{cases};\qquad \kappa(D_w)=\ck_{w^0}.
\end{equation}

The equality for $\kappa(T_w)$ follows from \cite[Proposition 2.18, Lemma 2.19]{ISY}. For the second equality, note that $\kappa(bT_i)=0$ for $1\le i\le n$ and $\kappa(bD_i)=\kappa(b+bT_i)=\kappa(b)$.

The Demazure product on $\W$ is defined by $w\star s_i=ws_i$ if
$\ell(ws_i)>\ell(w)$ and $w\star s_i=w$ otherwise, and extends to
$w\star v$ along any reduced expression for $v$. Using \cite[Lemma 2.10]{ISY}, we have
\begin{equation}\label{eq:demazure}
 T_wT_v=(-1)^{\ell(w)+\ell(v)-\ell(w\star v)}T_{w\star v},\quad D_wD_v=D_{w\star v}.
\end{equation}

\begin{lemma}\label{lem:TD}
For $z\in \W$, the relations between $T_z$ and $D_z$ are given by:
\begin{equation}\label{eq:TD}
 D_z=\sum_{y\le z}T_y,\qquad
 T_z=\sum_{y\le z}(-1)^{\ell(z)-\ell(y)}D_y.
\end{equation}
\end{lemma}
\begin{proof}
The first equality is proved in \cite[Proposition 2.11]{ISY}.
  For the second, define an involution on $\mathbb K_0$ by
\[
\omega(T_i)=-D_i=-1-T_i.
\]
 Since \(D_i^2=D_i\), then
$(\omega(T_i))^2=(-D_i)^2=D_i==-\omega(T_i)$,
so $\omega(T_i)$ satisfy the same quadratic relation and braid relations as \(T_i\).
Moreover,
\[
\omega(D_i)=1+\omega(T_i)=1-D_i=-T_i,
\]
and therefore \(\omega^2(T_i)=T_i\).
Now apply \(\omega\) to the first equality of \eqref{eq:TD} gives
\[
(-1)^{\ell(z)}T_z
=
\sum_{y\le z}(-1)^{\ell(y)}D_y.
\]
Multiplying both sides by \((-1)^{\ell(z)}\) proves the second equality.
\end{proof}

Applying $\kappa$ to the identities in \eqref{eq:TD} and using \eqref{eq:kappa}
gives 
\[
 \ell_w=\sum_{\substack{u\in W_{\mathrm{af}}^0\\u\le w}}k_u, \quad k_w=\sum_{\substack{u\in W_{\mathrm{af}}\\u\le w}}(-1)^{\ell(w)-\ell(u)}\ell_{u^0}
 \text{ for all }w\in W_{\mathrm{af}}^0.
\]
The first equality appears in \cite[Equation~(2.39)]{ISY}. Both the Pieri classes are related by
\[
 \cg_r=\sum_{a=0}^r\gh_a, \quad \gh_r=\cg_r-\cg_{r-1} \text{ for } 1\le r\le n.
\]

\subsection{Interpolation formulas for ordinary and closed Pieri classes}
We derive explicit formulas of both ordinary and closed Pieri classes in terms of translation elements using Lagrange interpolation. 

We recall $\theta^\vee=\eps_1-\eps_N$. For $1\le j\le n$, define $\tau_j:=t_{\eps_j-\eps_N}$. 

\begin{lemma}
 For $1\le j \le n$, a reduced expression for $\tau_j$ is
   \begin{equation}\label{eq:tauword}
 \tau_j=(s_{j-1}\cdots s_1)s_0(s_1\cdots s_n)(s_{n-1}\cdots s_j).
\end{equation} 
\end{lemma}
\begin{proof}
 Using
$s_\theta=s_1\cdots s_n s_{n-1}\cdots s_1$ in
$t_{\theta^\vee}=s_0s_\theta$, we have
\[
\tau_1= t_{\theta^\vee}=s_0s_1\cdots s_n s_{n-1}\cdots s_1.
\]

Conjugate this translation by
$p_j=s_{j-1}\cdots s_1$. Since $p_j(\eps_1-\eps_N)=\eps_j-\eps_N$,
$$\tau_j=t_{\eps_j-\eps_N}=p_jt_{\eps_1-\eps_N}p_j^{-1}.$$
The final $s_{j-1}\cdots s_1$ cancels with
$p_j^{-1}$ and the resulting word becomes \eqref{eq:tauword}.
\end{proof}

\begin{proposition}\label{prop:newton}
For $1\le j\le n$,
\begin{equation}\label{eq:newton}
 \frac{\tau_j-1}{1-y_j}
 =\sum_{a=1}^j\left(\prod_{i=1}^{a-1}(1-y_j/y_i)\right)\gh_a.
\end{equation}
\end{proposition}
\begin{proof}
Since $s_0\ast 1=\tau_1$ and $e^{\alpha_0}=e^{-\theta}=y_1$, we have
\[
 k_1=T_0\ast 1=\frac{\tau_1-1}{1-y_1}.
\]
This proves \eqref{eq:newton} for $j=1$ and completes the proof for $\mathrm{A}_1$. For $j<n$, $s_j$ interchanges
$y_j, y_{j+1}$, fixes the other $y_i$, and conjugates $\tau_j$ to
$\tau_{j+1}$. Also using $s_j=1+(1-e^{\alpha_j})T_j$
\[
 s_j\ast\gh_j=\gh_j+(1-y_{j+1}/y_j)\gh_{j+1},
\] as
$e^{\alpha_j}=y_{j+1}/y_j$. If $a<j$, 
we have
$s_j\ast\gh_a=\gh_a$ as $T_j\ast \gh_a=0$, see \cite[Proposition 2.24]{ISY}.

For $1\le j<n$, assume \eqref{eq:newton} for $j$ and apply $s_j\ast$. 
Then
\begin{align*}
 \frac{\tau_{j+1}-1}{1-y_{j+1}}
 ={}&\sum_{a=1}^{j}
       \left(\prod_{i=1}^{a-1}(1-y_{j+1}/y_i)\right)\gh_a+\left(\prod_{i=1}^{j-1}(1-y_{j+1}/y_i)\right)
       (1-y_{j+1}/y_j)\gh_{j+1}.
\end{align*}
This is precisely the asserted formula for $j+1$.
\end{proof}

\begin{proposition}
\label{thm:identify}
For every $1\le r\le n$, the Pieri classes are given by
\begin{equation}\label{eq:translation}
 \gh_r=(-1)^{r-1}Y_{r-1}
 \sum_{j=1}^r\frac{\tau_j-1}
 {(1-y_j)\prod_{k\le r,\,k\ne j}(y_j-y_k)}.
\end{equation}
\begin{equation}\label{eq:closedtranslation}
 \cg_r=1+(-1)^{r-1}Y_r\sum_{j=1}^r
 \frac{\tau_j-1}{y_j(1-y_j)\prod_{k\le r,\ k\ne j}(y_j-y_k)}.
\end{equation}
\end{proposition}

\begin{proof}
We work in $\mathbb L^\Delta
=R(T)^\Delta\otimes_{R(T)}\mathbb L$ 
and define the polynomial
\begin{equation}\label{eq_H_r}
     H_r(t)=\sum_{a=1}^{r}
 \left(\prod_{i=1}^{a-1}(1-t/y_i)\right)\gh_a\in\mathbb L^\Delta[t].
\end{equation}
$H_r(t)$ is a polynomial in $t$ of degree $r-1$. At $t=y_j$, all terms with $a>j$
vanish, and Proposition~\ref{prop:newton} gives
$H_r(y_j)=(\tau_j-1)/(1-y_j)$. Hence Lagrange interpolation gives
\begin{equation}\label{eq_H_r2}
     H_r(t)=\sum_{j=1}^{r}\frac{\tau_j-1}{1-y_j}
 \frac{\prod_{k\le r,\ k\ne j}(t-y_k)}
      {\prod_{k\le r,\ k\ne j}(y_j-y_k)}.
\end{equation}
Comparing the coefficients of $t^{r-1}$ in $H_r(t)$ from \eqref{eq_H_r} and \eqref{eq_H_r2} gives
\[
 \frac{(-1)^{r-1}}{Y_{r-1}}\gh_r
 =\sum_{j=1}^r
 \frac{\tau_j-1}{(1-y_j)\prod_{k\le r,\ k\ne j}(y_j-y_k)}.
\]
This proves \eqref{eq:translation}.

For the closed Pieri class using \eqref{eq_H_r}, $H_r(0)=\sum_{a=1}^r\gh_a=\cg_r-1$, while
\[
 \prod_{k\le r,\ k\ne j}(-y_k)=(-1)^{r-1}Y_r/y_j.
\]
This proves \eqref{eq:closedtranslation} using \eqref{eq_H_r2}.
\end{proof}

\begin{definition}
Let $\mathcal Z_A$ be the set of affine Weyl group elements those lie
below at least one of the translations $\tau_1,\ldots,\tau_n$ in Bruhat order.
Thus
\[
  \mathcal Z_A
  := \{w\in\W \mid w\le \tau_j \text{ for some } j\in\{1,\ldots,n\}\}\subset \W.
\]
\end{definition}

Then $\mathcal Z_A$ is a finite set as every element of $\mathcal Z_A$ have lengths at most $2n$.
Using \cite[Lemma 2.3, Proposition 2.4]{LSS} every $v\in \W\subset \mathbb K$ can be expressed as an $R(T)$-linear combination of $T_w$ for $w\leq v$.
Thus \eqref{eq:translation} and \eqref{eq:closedtranslation} can be written as $$\gh_r=\sum_{w\in \mathcal{Z}_A} C_{r,w} T_w, \quad \cg_r=\sum_{w\in \mathcal{Z}_A} E_{r,w} D_w=\sum_{w\in \mathcal{Z}_A} \CC_{r,w} T_w.$$
Now using Lemma \ref{lem:TD},
\begin{equation}\label{eq:E}
 E_{r,w}=\sum_{\substack{z\in \mathcal{Z}_A\\z\ge w}}(-1)^{\ell(z)-\ell(w)}\CC_{r,z},
 \qquad \CC_{r,z}=\sum_{a=0}^r C_{a,z},\quad C_{0,z}=\delta_{z,\id}.
\end{equation}

For $F\in\mathbb K$, let $[T_w]F$ and $[D_w]F$ denote the coefficients of $T_w$ and $D_w$ in the respective basis expansions of $F$.
We now compute $[T_w]\tau_j$ and $[D_w]\tau_j$ explicitly.

\section{An auxiliary algebra for translation coefficients}\label{sec_aux_alg}
Let $H_0$ be the associative unital $\mathbb Z$-algebra generated
by $\{U_i: i\in I\}$, with $U_i^2=-U_i$ and the Coxeter braid relations.
It has basis $\{U_w\}_{w\in \W}$, where $U_w$ is defined by \eqref{eq_pro} with $T_i$ is replaced by $U_i$. Define the auxiliary algebra
\[
 \mathcal H:=R(T)^\Delta\otimes_{\mathbb Z}H_0,
\]
where the product is defined by $(f\otimes u)(g\otimes v)=(fg\otimes uv)$.
We define a map
$$\varphi: \mathcal H\to \mathbb K^\Delta;\quad  \varphi(\sum_w f_wU_w)=\sum_w f_wT_w.$$ Then $\varphi$ is an $R(T)^\Delta$-linear identification
onto the localized $K$-nilHecke algebra. For $z\in R(T)^\Delta$, define $R_i(z):=1+(1-z)U_i$.

\begin{lemma}\label{lem:roots}\cite[Proposition 2.10]{LSS}
Let $v=s_{i_1}\cdots s_{i_m}\in \W$ be a reduced expression and
\[
\beta_1=\alpha_{i_1},~ \beta_k=s_{i_1}\cdots s_{i_{k-1}}(\alpha_{i_k}) \text{ for } 2\le k\le m.
\]
Restrict these affine roots to the finite torus when forming characters.
Then
\begin{equation}\label{eq:rootproduct}
 [T_w]v=[U_w]\prod_{k=1}^m R_{i_k}(e^{\beta_k}).
\end{equation}
For every $J=\{j_1<\cdots<j_r\}\subseteq\{1,\dots,m\}$, let
$w_J=s_{i_{j_1}}\star\cdots\star s_{i_{j_r}}$ be the Demazure product of the
corresponding subword, with $w_\varnothing=1$. Then
\begin{equation}\label{eq_billey}
 [T_w]v
 =\sum_{\substack{J\subseteq\{1,\dots,m\}\\ w_J=w}}
   (-1)^{|J|-\ell(w)}
   \prod_{k\in J}(1-e^{\beta_k}).
\end{equation}
\end{lemma}
\begin{proof}
Let $\sigma_i$ act on every element $\sum_w f_wU_w\in \mathcal H$ by 
$$\sigma_i\left(\sum_w f_wU_w\right)=\sum_w s_i(f_w)U_w.$$
Then, for every $F_0\in\mathcal H$,
\begin{equation}\label{eq_phi}
  \varphi\bigl(R_i(e^{\alpha_i})\sigma_i(F_0)\bigr)=s_i\varphi(F_0).  
\end{equation}

To consider this, consider $F_0=fU_w$. The twisted multiplication of $\mathbb K^\Delta$ as in \eqref{eq_tw_act}
gives $s_i f=s_i(f)s_i$. Using $s_i=1+(1-e^{\alpha_i})T_i$, we obtain
\[
\begin{aligned}
 s_i\varphi(F_0)
 =s_i(f)s_iT_w=s_i(f)T_w+(1-e^{\alpha_i})s_i(f)T_iT_w=\varphi\!\left(R_i(e^{\alpha_i})\sigma_i(F_0)\right).
\end{aligned}
\]
Apply the identity in \eqref{eq_phi} inductively complete the proof of \eqref{eq:rootproduct}.
Each new simple reflection acts on all subsequent root characters, giving
exactly the roots $\beta_k$. Now expanding the product in \eqref{eq:rootproduct} and using the $0$-Hecke multiplication proves
\eqref{eq_billey}.
\end{proof}

Formula \eqref{eq_billey} is the $K$-theoretic analogue of the formula in \cite{AJS, Bil} for the restriction of a torus-equivariant cohomology class $[X_v]$ to a torus-fixed point $w$.

\subsection{Yang–Baxter identities and the translation root product}
We have the following identities on $R_i(z)$
\begin{align}
 R_i(u)R_i(v)&=R_i(uv),\label{eq:quad}\\
 R_i(u)R_{i+1}(uv)R_i(v)
 &=R_{i+1}(v)R_i(uv)R_{i+1}(u).\label{eq:YB}
\end{align}
For $i\in I$, $R_i(z)$ are standard
$0$-Hecke Yang--Baxter operators in \cite[Section 4]{LLT}, $R_i(z)=Y_i^{\bar\pi}(1,z)$ with $U_i=\bar\pi_i$, see also
\cite{FK1,FK2}.

\begin{lemma}\label{lem:rootlist}
For $\tau_j$, the ordered list of $e^{\beta_k}$ is
\begin{equation}\label{eq:rootlist}
 \left(\frac{y_j}{y_{j-1}},\ldots,\frac{y_j}{y_1};\ y_j;\
 y_1,\ldots,y_{j-1},y_{j+1},\ldots,y_n,y_j;\
 \frac{y_j}{y_n},\ldots,\frac{y_j}{y_{j+1}}\right).
\end{equation}
\end{lemma}
\begin{proof}
In the four blocks of \eqref{eq:tauword}, the restricted roots $\beta_k$ are
\[
\begin{array}{c|c}
\text{reflection and its range}&\text{restricted root}\\ \hline
s_i\text{ in the initial block},\ j-1\ge i\ge1& a_i-a_j\\
s_0& a_N-a_j\\
s_i\text{ in the increasing block},\ 1\le i<j& a_N- a_i\\
s_i\text{ in the increasing block},\ j\le i<n& a_N-a_{i+1}\\
s_n\text{ in the increasing block}& a_N-a_j\\
s_i\text{ in the final block},\ n-1\ge i\ge j& a_{i+1}-a_j
\end{array}
\]
Writing characters $e^{\beta_k}$ using $y_i$ proves the list \eqref{eq:rootlist}.
\end{proof}

Let $\mathcal R_j$ be the ordered auxiliary product associated with
$\tau_j$. Products are taken in the indicated order, and empty
products are $1$.
\begin{equation}\label{eq:translation-root-product}
\begin{aligned}
 \mathcal R_j
 &=\Big(\prod_{i=j-1}^1R_i(\frac{y_j}{y_i})\Big)R_0(y_j)
   \Big(\prod_{i=1}^{j-1}R_i(y_i)\Big)
   \Big(\prod_{i=j}^{n-1}R_i(y_{i+1})\Big)R_n(y_j)
   \Big(\prod_{i=n-1}^jR_i(\frac{y_j}{y_{i+1}})\Big).
\end{aligned}
\end{equation}

Let us denote $B_j(w)=[T_w](\tau_j-1)$. Then for $w\neq id$, using Lemma \ref{lem:roots} 
$$B_j(w)=[U_w]\mathcal R_j.$$

In particular $B_j(\id)=0$. Let $\iota$ be the involution of $R(T)^\Delta$ defined by
$\iota(e^\lambda)=e^{-\lambda}$ for $\lambda\in X^*(T)$.
It extends coefficientwise to $\mathcal H$ by fixing every $U_w$.

\begin{lemma}\label{lem:jp-basis}
 For $w\ne\mathrm{id}$, we have
\begin{equation*}
 [D_w]\tau_j=(-1)^{\ell(w)}y_j^{n+1}\iota(B_j(w)),
 \qquad [D_{\mathrm{id}}]\tau_j=y_j^{n+1}.
\end{equation*}
\end{lemma}
\begin{proof}
In the auxiliary algebra $\mathcal H$ consider
$V_i=1+U_i$ and define $\omega(U_i)=-V_i$, fixing the coefficient ring $R(T)^\Delta$. This is an involution and $\omega(V_i)=-U_i$.
Direct computation gives
\[
 R_i(z)=z\,\omega(R_i(z^{-1})).
\]
Applying this identity in \eqref{eq:rootlist} gives
\[
 \mathcal R_j=y_j^{n+1}\omega(\iota(\mathcal R_j)).
\]
Since $\omega(U_w)=(-1)^{\ell(w)}V_w$, the $V_w$-coefficient of $\mathcal R_j$ is
$(-1)^{\ell(w)}y_j^{n+1}\iota([U_w]\mathcal R_j)$.
For $w\ne\mathrm{id}$ the last coefficient is $B_j(w)$,
for the identity it is $1$.

The coefficient wise linear identification $\varphi$
carries $V_w$ to $D_w$. Now the proof follows using $[V_w]\mathcal{R}_j=[D_w]\tau_j$.
\end{proof}

We define support of an element $w \in \mathcal Z_A$ by the set of indices appearing in a reduced expression of $w$ and denote this by $\Supp(w)$. Note that if $i\notin \Supp(w)$ then substituting $U_i=0$, (equivalently substituting $R_i(z)=1$) does not effect the coefficients of $U_w$ in $\mathcal R_j$. For a proper subset $S\subsetneq I$, we denote $\mathcal R_j|_S$ after setting $U_i=0$ in $\mathcal R_j$ for $i\notin S$.

\subsection{Auxiliary path products on cyclic intervals}
Let $C=[a,a+1,\ldots,b-1,b]$ be a nonempty proper cyclic interval of nodes. We denote this interval by $C=[a,b]$.
The variable $t$ is independent until explicitly evaluated. Define
\begin{equation}\label{eq:pathproduct}
 H_C(t):=R_a(y_{a+1})\cdots R_{b-1}(y_b)R_b(t)
             R_{b-1}(t/y_b)\cdots R_a(t/y_{a+1}).
\end{equation}
If $C=\{a\}$ is a singleton set then $H_{\{a\}}(t):=R_a(t)$.

\begin{lemma}\label{lem:products}
For a proper cyclic interval $[a,b]$, we have
\begin{align}
 H_{[a,b]}(t)
 &=R_b(t/y_b)\cdots R_{a+1}(t/y_{a+1})R_a(t)
   R_{a+1}(y_{a+1})\cdots R_b(y_b),
 \label{eq:reversepivot}\\
 H_{[a,b]}(y_j)
 &=H_{[a,j-1]}(y_j)H_{[j,b]}(y_j)
 \qquad(a<j\le b),
 \label{eq:cut}
\end{align}
where $a<j\le b$ refers to the cyclic order. Moreover, for every
proper subset $S\subsetneq I$,
\begin{equation}\label{eq:restrictedH}
 \left.\mathcal R_j\right|_S
 =\prod_{C:\text{ component of }S}\left.\mathcal R_j\right|_C
 =\prod_{C:\text{ component of }S}H_C(y_j).
\end{equation}
\end{lemma}
\begin{proof}
Applying \eqref{eq:YB} to the three central factors in
\eqref{eq:pathproduct}, and commuting the nonadjacent
factors gives
\[
 H_{[a,b]}(t)=R_b(t/y_b)H_{[a,b-1]}(t)R_b(y_b).
\]
Iteration proves \eqref{eq:reversepivot}. At $t=y_j$ the factor
$R_{j-1}(t/y_j)$ in \eqref{eq:pathproduct} becomes $1$,
so the factors on $[a,j-2]$ commute past those on $[j,b]$, giving \eqref{eq:cut}.

Setting $U_i=0$ for $i\notin S$ replaces $R_i(z)$ by $1$. Since factors
from distinct components commute,
this gives the first equality in
\eqref{eq:restrictedH}.

For a component $C=[a,b]$, that does not contain $0$, restrict \eqref{eq:rootlist}
to $C$. If $j\le a$, combining the two endpoint factors by
\eqref{eq:quad} gives \eqref{eq:pathproduct} at $t=y_j$.
For $j>b$ and $a<j\le b$, the restricted expression is given by
\eqref{eq:reversepivot} and \eqref{eq:cut}, respectively.
Thus $\left.\mathcal R_j\right|_C=H_C(y_j)$.

For a component $C=[a,\ldots,n,0,\ldots,b]$, if $j$ is internal, then apply
\eqref{eq:cut} to the interval $[a,n]$ or $[0,b]$ that contains $j$.
Since $y_0=1$, the factor $R_n(y_0)$ is trivial, we may use \eqref{eq:reversepivot} on
$[0,b]$ and commute this block past the increasing factors on $[a,n-1]$.
The remaining factors are $R_n(t)$ followed by the decreasing factors
on $[a,n-1]$. At $t=y_j$, this recovers the restriction of
\eqref{eq:rootlist}, proving $\left.\mathcal R_j\right|_C=H_C(y_j)$.
\end{proof}

\section{A characterization of \texorpdfstring{$\mathcal Z_A$}{\mathcal Z_A} and a formula for \texorpdfstring{$B_j(w)$}{B_j(w)}}\label{sec:colored}
 In this section, we introduce a combinatorial data on $S=\Supp(w)$ of each element $w\in \mathcal Z_A$ to characterize the elements of $\mathcal{Z}_A$. Moreover, we use this combinatorial data to compute the coefficient $B_j(w)$ for all $w\in \mathcal Z_A$.

Either $S=I$ for full support elements, or $S$ is proper, its connected components are cyclic intervals, with indices taken modulo $n+1$. 
 We call an edge joining $i-1$ and $i$ an internal edge of $S$ if both $i-1$ and $i$ lies in $S$. 
 We call $i$ as the head of the edge joining $i$ and $i-1$ and index this edge by its head $i$. The edge between $n$ and $0$ can be also considered in this way and index this edge by its head $0$. On each heads of internal edges of $S$ we introduce label $L$, $R$, or, both of $L$ and $R$.  we also label those edges by $B$, if it gets both the label $L$ and $R$.

\subsection{Support data for proper support}

Let $C=[a,a+1,\ldots,b]$ be a nonempty proper cyclic interval, and $E(C)=C\setminus\{a\}$ be its set of internal edge heads.
For $L_C,R_C\subseteq E(C)$ with $L_C\cup R_C=E(C)$, set
\begin{equation}\label{eq:path-word}
 w(L_C,R_C):=
 \left(
   \prod_{\substack{i\in C\setminus\{b\}\\i+1\in L_C}}^{\longrightarrow}
   s_i
 \right)
 s_b
 \left(
   \prod_{\substack{i\in C\setminus\{b\}\\i+1\in R_C}}^{\longleftarrow}
   s_i
 \right).
\end{equation}
The arrows indicate increasing and decreasing path order.
Thus each internal head has a left choice, a right choice, or both. In the terminology of \cite[Section~5.1]{LS}, the expression for
$w(L_C,R_C)$ in \eqref{eq:path-word} is a $\Lambda$-word. Distinct pairs $(L_C,R_C)$ give distinct elements $w(L_C,R_C)$ using \cite[Lemma~16]{LS}.

\begin{lemma}\label{lem:path-word}
The expression \eqref{eq:path-word} is reduced, has support $C$, and has
length $1+|L_C|+|R_C|$. 
Moreover, $w(L_C,R_C)$ uniquely determines $(L_C,R_C)$.
\end{lemma}

\begin{proof}
We identify $\langle s_i\mid i\in C\rangle$ with the symmetric group $S_{|C|+1}$.
Reducedness follows by induction on $|C|$.
The singleton case is immediate. Otherwise, let $p=(1,2)$ and
$v$ be the element given by the expression on the shorter path. By induction,
this expression is reduced and has full support on that path. Thus $v$ fixes
$1$ and both $v(2), v^{-1}(2)>2$.

The three choices in \eqref{eq:path-word} give $pv$, $vp$, and $pvp$.
Using \cite[(1.26), Proposition~1.5.2]{BB}
\[
 \ell(pv)=\ell(vp)=\ell(v)+1,\qquad
 \ell(pvp)=\ell(v)+2,
\]
where the last equality uses $(pv)(1)=2<(pv)(2)=v(2)$.
Hence the expression is reduced. Its support and length follow from
\eqref{eq:path-word}.

The three choices can be recovered from $w=w(L_C,R_C)$. If $w=pv$ then $w(1)=2$, if $w=vp$ then $w^{-1}(1)=2$. If $w=pvp$ then \(w(1)=v(2)>2\) and
\(w^{-1}(1)=v^{-1}(2)>2\). Now removing the corresponding
copies of $p$ recovers $v$, so induction determines the data $(L_C,R_C)$.
\end{proof}

\begin{lemma}\label{lem:path-coefficient}
Considering $w_C=w(L_C,R_C)$, $H_C(t)$ satisfies
\begin{equation}\label{eq:path-coefficient}
 [U_{w_C}]H_C(t)
 =
 (1-t)
 \prod_{i\in L_C}(1-y_i)
 \prod_{i\in R_C}(1-t/y_i).
\end{equation}
Every other coefficient indexed by an element of support exactly
$C$ is zero.
\end{lemma}

\begin{proof}
A selection contributing an element of support $C$ must use the
central $U_b$ and at least one occurrence of every other generator.
By Lemma~\ref{lem:path-word}, all such selections are reduced and
represent distinct elements. The central factor contributes $1-t$,
while the left and right choices at head $i$ contribute $1-y_i$ and
$1-t/y_i$, respectively. Multiplying these weights proves the formula.
\end{proof}

Let $\varnothing\ne S\subsetneq I$, and write
\[
 E(S):=\{i\in I:i-1,i\in S\}.
\]
An \emph{admissible support datum} is a triple $(S,L,R)$ satisfying
\[
 L,R\subseteq E(S),\qquad L\cup R=E(S),\qquad
 0\in E(S)\ \Longrightarrow\ 0\in R\setminus L.
\]
The condition at $0$ agrees with \eqref{eq:path-coefficient}, since
$y_0=1$ \eqref{eq:path-coefficient} vanish if $0\in L$.

For each component $C$ of $S$, put $L_C=L\cap E(C)$ and
$R_C=R\cap E(C)$, and define
\[
 w(S,L,R):=\prod_{C:\text{ component of }S}w(L_C,R_C).
\]
Distinct components commute. If $c$ is the number of components,
\[
 \operatorname{Supp}(w(S,L,R))=S,\qquad
 \ell(w(S,L,R))=c+|L|+|R|.
\]
Lemma~\ref{lem:path-word} proves that the admissible support datum $(S,L,R)$ is uniquely determined by its element $w(S,L,R)$.

\subsection{Support data for full support}
For full support, $S=I$. In this case, an \emph{admissible support datum} is a triple
$(I,L,R)$ satisfying
\[
 L\cup R=I,\qquad 0\in R\setminus L,\qquad L\setminus R\ne\varnothing.
\]
Thus there is at least one finite head labeled $L$.
Choose such a head $j\in L\setminus R$, call it a \emph{pivot}, and define
\begin{equation}\label{eq:fullword}
w(I,L,R):= \left(\prod_{\substack{i=j-1,\ldots,1\\i\in R}}s_i\right)s_0
 \left(\prod_{i=1}^{n-1}s_i^{\eta_i}\right)s_n
 \left(\prod_{\substack{i=n-1,\ldots,j\\i+1\in R}}s_i\right);~
 \eta_i=\begin{cases}\boldsymbol1_{i\in L},&i<j,\\
 \boldsymbol1_{i+1\in L},&i\ge j.
 \end{cases}
\end{equation}
The middle product is increasing. For the full support case, we set $c=0$. 

\begin{lemma}\label{lem:full}
Let $(I,L,R)$ be an admissible support datum.
\begin{enumerate}
\item The word \eqref{eq:fullword}
is reduced, has full support, and represents an element
$w(I,L,R)$ independent of the choice of pivot $j$. Its length is $|L|+|R|$.
\item Distinct admissible full-support data give distinct elements.
Thus $L,R$ are determined by $w(I,L,R)$.
\end{enumerate}
\end{lemma}

\begin{proof}
Let $w_j$ denote the word in~\eqref{eq:fullword} with pivot $j$. To prove pivot
independence, take consecutive pivots $a<b$ in $L\setminus R$ and set
$D=s_{b-1}\cdots s_a$ and $M=(L\cap R)\cap(a,b)$. After commuting
nonadjacent generators, $w_b$ and $w_a$ have the form
\[
  w_b=P D\Bigl(\prod_{i\in M}^{\longrightarrow}s_i\Bigr)Q,
  \qquad
  w_a=P\Bigl(\prod_{i\in M}^{\longrightarrow}s_{i-1}\Bigr)D Q,
\]
with the same outside elements $P,Q$, supported on $[0,a-1]$ and $[b,n]$ respectively. Now $w_a$ and $w_b$ agree since $Ds_i=s_{i-1}D$ for $a<i<b$.

Now take $a=\min(L\setminus R)$. Since $1,\ldots,a-1$ all belong to
$R$, the word is $w_a=\rho_a v$ with $v\in W$. The suffix $v$
contains $s_a,\ldots,s_n$. On the support component containing $n$,
it has the form~\eqref{eq:path-word}, and on each other component it is
increasing. Thus Lemma~\ref{lem:path-word} gives reducedness of $v$ and recovers its
increasing and decreasing selections. Since $\rho_a\in W_{\mathrm{af}}^0$,
the factorization $\rho_a v$ is length-additive. Thus $w_a$ is reduced and of full support and $\ell(w_a)=|L|+|R|$.

Finally, uniqueness of parabolic factorization recovers $\rho_a$ and
$v$ from $w_a$, hence also $a$. For $1\le i<a$, membership in $R$ is
forced, and $i\in L$ exactly when $s_i$ occurs in the increasing part
of $v$. For $a<i\le n$, the increasing and decreasing selections of
$s_{i-1}$ determine membership in $L$ and $R$, respectively. Together
with $a\in L\setminus R$ and $0\in R\setminus L$, this recovers the
entire datum.
\end{proof}

\subsection{Classification of \texorpdfstring{$\mathcal{Z}_A$}{\mathcal{Z}_A} using admissible support data} We classify the elements of $\mathcal Z_A$ in terms of admissible support data.
\begin{lemma}\label{lem:support-classification} 
The assignment $(S,L,R)\longmapsto w(S,L,R)$
is a bijection from the admissible support data with
$\varnothing\ne S\subseteq I$ to
$\mathcal Z_A\setminus\{\mathrm{id}\}$.
For the corresponding element $w=w(S,L,R)$,
\[
 \Supp(w)=S,\qquad
 \ell(w)=c+|L|+|R|,
\]
where $c$ is the number of components of $S$ when $S\subsetneq I$,
and $c=0$ when $S=I$.
\end{lemma}

\begin{proof}
The support and length assertions follow from
Lemma~\ref{lem:path-word} and Lemma~\ref{lem:full}.
These lemmas also give injectivity within each support case,
and the two cases have different supports.
It remains to identify the image with
$\mathcal Z_A\setminus\{\mathrm{id}\}$.

First suppose that $w\le\tau_j$ has nonempty proper support.
By \cite[Lemma~19]{LS}, it factors into connected $V$-words
on $[0,j-1]$ followed by connected $\Lambda$-words on $[j,n]$.
Convert the $V$-words into $\Lambda$-words using
\cite[Lemma~17]{LS}. Commutations merge adjacent factors
into the path form \eqref{eq:path-word}, in the cyclic order
of each support component. The standard word for $\tau_j$ places every $s_0$ before every $s_n$. This property passes to reduced subwords as $s_0s_n\neq s_ns_0$. Thus an internal head $0$ belongs to $R\setminus L$, and the
resulting datum is admissible.

Conversely, let $(S,L,R)$ be admissible with $S\subsetneq I$.
Choose $j\in\{1,\ldots,n\}\setminus S$ if possible, and otherwise
choose $j=1$, so that $S=\{1,\ldots,n\}$.
For a component including 0, the condition $0\notin L$ allows the factors on $[0,j-1]$ to precede those on
$[j,n]$. Convert the former into $V$-words using
\cite[Lemma~17]{LS}, retaining $\Lambda$-words for the latter. 
Hence $w(S,L,R)\le\tau_j$.

Now suppose that $S=I$.
For any admissible datum, choose a pivot $j\in L\setminus R$.
The reduced word \eqref{eq:fullword} is a subword
of the displayed expression for $\tau_j$.
Thus $w(I,L,R)\in\mathcal Z_A$.

Conversely, suppose that $w\le\tau_j$ has full support.
Choose a reduced subword for $w$ in $\tau_j$ as in \eqref{eq:tauword}. It must select the unique occurrences of $s_0$
and $s_n$, and at least one occurrence of every other generator.
Assign head $0$ the label $R$ and head $j$ the label $L$.
For $1\le i<j$, assign head $i$ the label $R$ when $s_i$
is selected in the initial decreasing block, and the label $L$
when it is selected in the middle increasing block.
For $j<i\le n$, assign head $i$ the label $L$ when $s_{i-1}$
is selected in the middle block, and the label $R$ when it is
selected in the final decreasing block.

Every head receives at least one label. Moreover,
$0\in R\setminus L$ and $j\in L\setminus R$, so the datum
is admissible. The selected word is exactly
\eqref{eq:fullword} with pivot $j$.
Therefore $w=w(I,L,R)$, completing the proof.
\end{proof}

\subsection{A product formula for \texorpdfstring{\(B_j(w)\)}{B_j(w)}}

We now describe a formula of $B_j(w)$ using
Lemma~\ref{lem:products}, Lemma~\ref{lem:path-word} for proper support and
Lemma~\ref{lem:full} for full support.
For an admissible support datum $(S,L,R)$, set
\[
 m=\ell(w)=c+|L|+|R|,\qquad Z=R\setminus\{0\},\qquad
 h=c-1+\boldsymbol1_{0\in R}.
\]
Then $h\ge0$ as $c\ge1$ for proper support and $\mathbf1_{0\in R}=1$ for full support. We call $w$ is \emph{admissible} if $w=w(S,L,R)$ for some admissible
datum $(S,L,R)$. Otherwise we call $w$ is \emph{non-admissible}. Corresponding to every $w=w(S,L,R)$, define a polynomial
\begin{equation}\label{eq:p}
 p_w(t):=\prod_{i\in L}(1-y_i)(1-t)^h
               \prod_{i\in Z}(1-t/y_i).
\end{equation}
Lemma \ref{lem:path-word} and Lemma \ref{lem:full} imply that $w$ determines the data $(S,L,R)$ and hence $p_w(t)$.
The number of factors in $p_w(t)$ is $m-1$ (including the constant factors $(1-y_i)$) as
\[
 |L|+h+|Z|=c+|L|+|R|-1=m-1.
\]
\begin{lemma}\label{lem:colored}
For every admissible element $w=w(S,L,R)$ the coefficient $B_j(w)$ is
\begin{equation}\label{eq:Bfactor}
B_j(w)=(1-y_j)p_w(y_j)\qquad(1\le j\le n).
\end{equation}
For every other element, including the identity, $B_j(w)=0$ for all $j$.
\end{lemma}
\begin{proof}
 We compute $B_j(w)=[U_w]\mathcal R_j$ for all admissible $w\ne\id$. If $\Supp(w)\subsetneq I$, write $w=\prod_{C} w_C$ for disjoint components $C$. Lemma~\ref{lem:path-coefficient} and \eqref{eq:restrictedH} implies  
\begin{align*}
 B_j(w)
 &=\prod_C[U_{w_C}]H_C(y_j)\\
 &=\prod_C\left[(1-y_j)
       \prod_{i\in L_C}(1-y_i)\prod_{i\in R_C}(1-y_j/y_i)\right]\\
 &=(1-y_j)^c\prod_{i\in L}(1-y_i)
                  \prod_{i\in R}(1-y_j/y_i)\\
 &=(1-y_j)^{c+\boldsymbol1_{0\in R}}
       \prod_{i\in L}(1-y_i)\prod_{i\in Z}(1-y_j/y_i)
  =(1-y_j)p_w(y_j).
\end{align*}

If $w=w(I,L,R)$, fix $j\in \{1,2,\dots,n\}$.
If $j\notin R$ then using \eqref{eq:fullword} and Lemma~\ref{lem:full}, $w$ has
exactly one selection as a subword of $\tau_j$. The unique selection's forced
$s_0,s_n$ contribute $(1-y_j)^2$. By \eqref{eq:rootlist}, each remaining head in $L$ contributes $1-y_i$ and each remaining
head in $R$ contributes $1-y_j/y_i$, with both factors for a head
labeled $B$. Hence
\[
 B_j(w)=(1-y_j)^2
       \prod_{i\in L\setminus\{j\}}(1-y_i)
       \prod_{i\in Z}(1-y_j/y_i)
       =(1-y_j)p_w(y_j),
\]
because $j\in L$, $c=h=0$, and $0\in R$.

If $j\in R$ then using Lemma~\ref{lem:full}, $w$ cannot be a subword of $\tau_j$,
so $B_j(w)=0$; the factor $1-y_j/y_j$ makes the asserted expression
zero as well.

The empty selection contributes to $U_{\id}$, subtracting $1$
therefore gives $B_j(\id)=0$.
\qedhere
\end{proof}

\begin{corollary}\label{cor:translation-containment}
For $1\le j\le n$ and
$w\in W_{\mathrm{af}}\setminus\{\mathrm{id}\}$,
\begin{equation}\label{eq_1}
    B_j(w)\ne0 \quad\Longleftrightarrow\quad w\le\tau_j. 
\end{equation}

If $w=w(S,L,R)$ has admissible support data, then 
\[
 w\le\tau_j \quad\Longleftrightarrow\quad j\notin R.
\]
\end{corollary}

\begin{proof}
The equivalent condition in \eqref{eq_1} follows from definition of $B_j(w)$ and Lemma \ref{lem:roots}. 
Let  $w=w(S,L,R)$ has admissible support data, the coefficient formula gives
\[
 B_j(w)=(1-y_j)^{h+1}
 \prod_{i\in L}(1-y_i)
 \prod_{i\in R\setminus\{0\}}(1-y_j/y_i).
\]
If $j\in R$, the last product contains the zero factor
$1-y_j/y_j$. If $j\notin R$, every factor in $B_j(w)$ is nonzero, since $0\notin L$.
Therefore $B_j(w)\ne0$ exactly when $j\notin R$.
\end{proof}

\section{Explicit formulas of the coefficients of ordinary and closed Pieri classes}\label{sec:marks}

In this section, we derive explicit formulas for
$C_{r,w}$, $\widetilde C_{r,w}$, and $E_{r,w}$.
The proofs use a common finite-sum identity, together with
the translation formulas and admissible support data
established in Section \ref{sec:colored}.

Let $w\ne\id$ have admissible support data $(S,L,R)$ such that $w=w(S,L,R)$. Write
\begin{equation}\label{eq:jp-data}
 \begin{gathered}
 a=n-c-|R|, \qquad A=Z\cap[1,r],\qquad
 J=[1,r]\setminus A=\{j_1<\cdots<j_k\}.
 \end{gathered}
\end{equation}

For proper support, $|R|\le |S|-c\le n-c$ and for full support,
$|R|\le n$ since $L\setminus R\ne\varnothing$. These inequalities
give $a\ge0$.

If $Z\setminus A=\{i_1<\cdots<i_s\}$, recall the ordered list $(b_1,\dots,b_d)$ defined in \eqref{eq:jp-list}.
Then
\begin{equation}\label{eq:jp-counts}
 \begin{gathered}
 k=r-|A|,\qquad
 d=h+|Z\setminus A|=m-1-|L|-|A|,\\
 a+h+|Z|=n-1,\qquad a+d=n-1-|A|.
 \end{gathered}
\end{equation}
Each $b_t$ is either $0$ or greater than $r$. 
For a subset $M$ of positive integers, put
\begin{equation}\label{eq:jp-M}
 M_\ell=1+|M\cap[1,\ell-1]|.
\end{equation}
If $|M|=k-1$, then $1\le M_\ell\le k$.
For example, if $k=3$, $d=5$, and $M=\{2,4\}$, then 
\[
(M_1,M_2,M_3,M_4,M_5)=(1,1,2,2,3).
\]

We use $[1,0]=\varnothing$, empty products equal to $1$, and empty
sums equal to $0$. The following lemma is the specialization of the general Leibniz formula
for divided differences~\cite[Corollary~30]{deBoor}.

\begin{lemma}\label{lem:jp-finite}
Let $k\ge1$, $p\ge0$, and let $z_1,\ldots,z_k$ lie in a commutative
ring $\mathscr R$ with identity, with $z_v-z_u$ a unit for $u\ne v$.
For $f_\ell(t)=u_\ell t+v_\ell\in\mathscr R[t]$, we have
\begin{equation}\label{eq:jp-finite}
 \sum_{v=1}^k
 \frac{\prod_{\ell=1}^p f_\ell(z_v)}
      {\prod_{u\ne v}(z_v-z_u)}
 =\sum_{\substack{M\subseteq[1,p]\\|M|=k-1}}
    \left(\prod_{\ell\in M}u_\ell\right)
    \left(\prod_{\ell\notin M}f_\ell(z_{M_\ell})\right).
\end{equation}
In particular, the expression vanishes when $p<k-1$.
\end{lemma}
\begin{proof}
Let $F_{k,p}$ and $G_{k,p}$ denote the left- and right-hand sides of \eqref{eq:jp-finite}, respectively, formed from the first $k$ nodes and the first $p$ factors. 

For $k=1$, only $M=\varnothing$ contributes, so both sides equal $\prod_{\ell=1}^p f_\ell(z_1)$. For $p=0$ and $k>1$, the empty set has no subset of size $k-1$, so $G_{k,0}=0$. Also, $F_{k,0}=0$ by comparing the coefficients of $t^{k-1}$ in the Lagrange interpolation identity
\[
1=\sum_{v=1}^k
\frac{\prod_{u\ne v}(t-z_u)}{\prod_{u\ne v}(z_v-z_u)}.
\]

Now suppose $k\ge2$ and $p\ge1$. Substituting
\[
f_p(z_v)=f_p(z_k)+u_p(z_v-z_k)
\]
into $F_{k,p}$ gives
\begin{equation}\label{eq:jp-recursion}
F_{k,p}=f_p(z_k)F_{k,p-1}+u_pF_{k-1,p-1}.
\end{equation}
Indeed, the constant term $f_p(z_k)$ gives the first summand. In the remaining sum, the term with $v=k$ vanishes, and for $v<k$ the factor $z_v-z_k$ cancels from the denominator, leaving the sum over the first $k-1$ nodes.

For $G_{k,p}$, split the sum according to whether $p\in M$. If $p\notin M$, then $M_p=k$ and the last factor is $f_p(z_k)$. These sets contribute $f_p(z_k)G_{k,p-1}$. If $p\in M$, the last factor is $u_p$. The subset $M\setminus \{p\}$ has size $k-2$ and removing $\{p\}$ preserve $M_\ell$ for every $\ell<p$. These sets therefore contribute $u_pG_{k-1,p-1}$.

The initial values and recurrences agree, so induction on $p$ proves \eqref{eq:jp-finite}. When $p<k-1$, $[1,p]$ has no subset of size $k-1$, so the both sides vanish.
\end{proof}

\begin{theorem}\label{thm:jp-ordinary}
Let $ w \ne \mathrm{id}$ be an element corresponding to admissible triple $(S,L,R)$.
If $k=0$ or, $d<k-1$ then $C_{r,w}=0$. Otherwise
\begin{equation}\label{eq:jp-ordinary}
 {\begin{aligned}
 (-1)^{\ell(w)-r}C_{r,w}
 =\frac{Y_{r-1}\prod_{i\in L}(y_i-1)}{\prod_{i\in A}y_i}\sum_{\substack{M\subseteq[1,d]\\|M|=k-1}} \left(\prod_{\ell\in M}\frac{1}{y_{b_\ell}}\right)
  \left(\prod_{\substack{1\le\ell\le d\\\ell\notin M}}
      \frac{y_{j_{M_\ell}}}{y_{b_\ell}}-1\right).
 \end{aligned}}
\end{equation}
Each complete summand, including its prefactor on the right hand side
has the form
\[
 x_1^{u_1}\cdots x_n^{u_n}
 \prod_{v=1}^{\ell(w)-r}(e^{-\beta_v}-1),
 \qquad u_i\in\mathbb Z_{\ge0},\quad \beta_v\in\Phi^+.
\]
If $w=\id$ or, $w$ is non-admissible then $C_{r,w}=0$. Consequently,
$(-1)^{\ell(w)-r}C_{r,w}\in\mathcal P$ for every $w\in \W$.
\end{theorem}
\begin{proof}
Taking the coefficient of $T_w$ in \eqref{eq:translation} and using
$B_j(w)=(1-y_j)p_w(y_j)$ from Lemma~\ref{lem:colored} gives
\[
 C_{r,w}=(-1)^{r-1}Y_{r-1}
 \sum_{j=1}^r
 \frac{p_w(y_j)}{\prod_{v\le r,\ v\ne j}(y_j-y_v)}.
\]
Since $p_w(y_j)=0$ for $j\in A$, the sum vanishes when $k=0$.
For $j\in J$, cancel each factor indexed by $i\in A$ using
\[
 \frac{1-y_j/y_i}{y_j-y_i}=-\frac1{y_i}.
\]
Change each remaining factor $1-y_i$ to $-(y_i-1)$ and each
$1-y_j/y_{b_\ell}$ to $-(y_j/y_{b_\ell}-1)$.
The total sign is
\[
 (-1)^{r-1+|A|+|L|+d}=(-1)^{m-r},
\]
by \eqref{eq:jp-counts}. Therefore
\begin{equation}\label{eq:jp-ordinary-sum}
 (-1)^{m-r}C_{r,w}
 =\frac{Y_{r-1}\prod_{i\in L}(y_i-1)}{\prod_{i\in A}y_i}
 \sum_{v=1}^k
 \frac{\prod_{\ell=1}^d(y_{j_v}/y_{b_\ell}-1)}
      {\prod_{u\ne v}(y_{j_v}-y_{j_u})}.
\end{equation}
Apply Lemma~\ref{lem:jp-finite} with $z_v=y_{j_v}$ and
$f_\ell(t)=t/y_{b_\ell}-1$, with $u_\ell=y_{b_\ell}^{-1}$.
This proves \eqref{eq:jp-ordinary}. The sum is empty if $d<k-1$.
The other vanishings follow from Lemma \ref{lem:colored}.

It remains to describe a nonzero summand. Since $j_{M_\ell}\le r$ and
each $b_\ell$ is either $0$ or greater than $r$, each ratio
$y_{j_{M_\ell}}/y_{b_\ell}$ is $e^{-\beta}$ for some $\beta\in\Phi^+$.
The same holds for each $y_i$ with $i\in L$ in a nonzero summand.
Thus there are $|L|+d-(k-1)=m-r$ factors of the form $ e^{-\beta}-1$.

After omitting $y_0=1$, the monomial in the denominator has distinct indices, denote these by
$e_1<\cdots<e_t$, with $t\le |A|+k-1=r-1$. Since $e_v\ge v$,
\[
 \frac{Y_{r-1}}{y_{e_1}\cdots y_{e_t}}
 =\prod_{v=1}^t\frac{y_v}{y_{e_v}}
  \prod_{v=t+1}^{r-1}y_v
\]
is a monomial in the $x_i$ with nonnegative exponents. This proves
the asserted form. Finally, $x_i=1+(x_i-1)$ and, for
$\beta=\alpha_p+\cdots+\alpha_q$,
\[
 e^{-\beta}-1=\prod_{i=p}^q\bigl(1+(x_i-1)\bigr)-1\in\mathcal P.
\]
Hence every summand belongs to $\mathcal P$.
\end{proof}

We recall the augmentation map $\phi_0:R(T)\to\mathbb Z$ defined by $\phi_0(e^\lambda)=1$.

\begin{corollary}\label{cor:jp-length}
Let $m=\ell(w)$ and let $\mathfrak m
 =\ker(\phi_0)$
be the augmentation ideal. Then $C_{r,w}=0$ if $m<r$
or, $m>2n$. For $m\ge r$,
\[
 C_{r,w}\in\mathfrak m^{\,m-r},
 \qquad \mathfrak m^0=R(T).
\]
In particular, $\phi_0(C_{r,w})=0$ unless $m=r$.
\end{corollary}

\begin{proof}
The sum in \eqref{eq:jp-ordinary} is empty unless $d\ge k-1$.
By \eqref{eq:jp-counts}, this implies $m-r\ge |L|\ge0$.
Each complete summand in Theorem~\ref{thm:jp-ordinary} contains exactly
$m-r$ factors of the form $e^{-\beta}-1$.
Since these factors belong to $\mathfrak m$, we obtain
$C_{r,w}\in\mathfrak m^{\,m-r}$.
Coefficients outside the admissible nonidentity elements
already vanish.

Each translation $\tau_j$ has length $2n$, so its
root-subword expansion contains no element of length
greater than $2n$. This proves the upper bound.
Applying $\phi_0$ proves the final assertion.
\end{proof}

\begin{theorem}\label{thm:closedinterp}
For a nonidentity admissible element $w$,
\begin{equation}\label{eq:closedfinite}
 \CC_{r,w}=(-1)^{r-1}Y_r\sum_{j=1}^r
 \frac{p_w(y_j)}{y_j\prod_{k\le r,\ k\ne j}(y_j-y_k)}.
\end{equation}
If $p_w(t)=\sum_{b=0}^{D}p_{w,b}t^b$, then
\begin{equation}\label{eq:closedpolynomial}
 \CC_{r,w}=p_w(0)+(-1)^{r-1}Y_r
 \sum_{b=r}^{D}p_{w,b}\mathsf h_{b-r}(y_1,\ldots,y_r).
\end{equation}
Also $\CC_{r,\id}=1$, and the other coefficients vanish.
\end{theorem}
\begin{proof}
Taking the coefficient of $T_w$ in
\eqref{eq:closedtranslation} and using
$B_j(w)=(1-y_j)p_w(y_j)$ from Lemma~\ref{lem:colored}
proves \eqref{eq:closedfinite}. The identity coefficient is $1$,
and the other vanishings follow from the same lemma.
For distinct $z_1,\ldots,z_k$ the elementary identities are

\begin{align}
 \sum_a\frac{z_a^{-1}}{\prod_{b\ne a}(z_a-z_b)}
 &=\frac{(-1)^{k-1}}{\prod_a z_a},\label{eq:interpnegative}\\
 \sum_a\frac{z_a^{v-1}}{\prod_{b\ne a}(z_a-z_b)}
 &=\begin{cases}0,&1\le v<k,\\ \mathsf h_{v-k}(z),&v\ge k.
 \end{cases}\label{eq:interppowers}
\end{align}
The first is Lagrange interpolation of the constant $1$ evaluated at
$0$. For the second compare coefficients in
\[
 \sum_a\frac{1}{\prod_{b\ne a}(z_a-z_b)(1-z_at)}
 =\frac{t^{k-1}}{\prod_a(1-z_at)};
\]
the rational identity follows by comparing its simple-pole terms.
Expanding $p_w$ in \eqref{eq:closedfinite} and applying these two
identities proves \eqref{eq:closedpolynomial}.
\end{proof}

\begin{corollary}
    If $p_w(t)$ is of degree less than $r$ then $\CC_{r,w}=p_w(0).$ 
\end{corollary}

\begin{remark}
Using Theorem~\ref{thm:jp-ordinary}, the coefficients $(-1)^{\ell(w)-r}C_{r,w}$ have exactly $\ell(w)-r$ factors $e^{-\beta}-1$. Considering $\mathcal P_-=\mathbb Z_{\ge0}[x_1,\ldots,x_n,1-x_1,\ldots,1-x_n]$, gives $C_{r,w}\in \mathcal P_-$, as replacing each factor $e^{-\beta}-1$ by  $-(1-e^{-\beta})$ cancels the sign.
Then using \eqref{eq:E} we obtain $\CC_{r,w}\in \mathcal P_-$. 
But $(-1)^{\ell(w)-r}\CC_{r,w}$ need not be in $\mathcal P$.

In Example \ref{eg_A2} for $\mathrm{A}_2$, $[T_{s_0s_1}]\cg_2=1-x_1x_2$. Here $\ell(01)=2$ and $1-x_1x_2\notin \mathcal P$ but $1-x_1x_2=(1-x_1)+x_1(1-x_2)\in \mathcal P_-$. 
However, $[D_{s_0s_1}]\ell_2= x_2(x_1x_2-1)\in\mathcal P$.
\end{remark}

\begin{theorem}
\label{thm:jp-closed}
Let $w\ne\mathrm{id}$ be an admissible element.
If $k=0$, then $E_{r,w}=0$. If $k\ge1$, then
\begin{equation}\label{eq:jp-closed}
 {\begin{aligned}
 (-1)^{m-r}E_{r,w}
 ={}&\frac{Y_r\prod_{i\in L}(y_i-1)}{\prod_{i\in L}y_i}\sum_{\substack{M\subseteq[1,a+d]\\|M|=k-1}}
 \left(\prod_{\substack{1\le\ell\le a\\\ell\notin M}}
 y_{j_{M_\ell}}\right)
 \left(\prod_{\substack{1\le t\le d\\a+t\notin M}}
 (y_{j_{M_{a+t}}}-y_{b_t})\right).
 \end{aligned}}
\end{equation}

If $w=\mathrm{id}$ then
\begin{equation}\label{eq:jp-closed-id}
 (-1)^rE_{r,\mathrm{id}}
 =Y_r\sum_{s=0}^{n-r}\mathsf h_s(y_1,\ldots,y_r)\in\mathcal P.
\end{equation}
If $w$ is non-admissible then $E_{r,w}=0$. In particular, $(-1)^{\ell(w)-r}E_{r,w}\in\mathcal P$
for every $w\in \W$.
\end{theorem}
\begin{proof}
Define
\[
 p_w^\vee(t):=\prod_{i\in L}(1-y_i^{-1})(1-t)^h\prod_{i\in Z}(1-ty_i).
\]
Then $\iota(B_j(w))=(1-y_j^{-1})p_w^\vee(y_j^{-1})$.
We apply Lemma~\ref{lem:jp-basis} to \eqref{eq:closedtranslation}. Then
\begin{equation}\label{eq:jp-closed-sum1}
 (-1)^{m-r}E_{r,w}
 =Y_r\sum_{j=1}^r
 \frac{y_j^{n-1}p_w^\vee(y_j^{-1})}
      {\prod_{v\le r,\ v\ne j}(y_j-y_v)}.
\end{equation}
For each summand,
\[
 y_j^{n-1}p_w^\vee(y_j^{-1})
 =\frac{\prod_{i\in L}(y_i-1)}{\prod_{i\in L}y_i}
   y_j^a(y_j-1)^h\prod_{i\in Z}(y_j-y_i),
\]
since $n-1-h-|Z|=a$.
Nodes in $A$ give zero. At a node in $J$, the factors $y_j-y_i$
for $i\in A$ cancel from numerator and denominator. Hence, for $k\ge1$,
\begin{equation}\label{eq:jp-closed-sum2}
 \begin{split}
 (-1)^{m-r}E_{r,w}
 ={}&\frac{Y_r\prod_{i\in L}(y_i-1)}{\prod_{i\in L}y_i}\sum_{v=1}^k
 \frac{y_{j_v}^{a}\prod_{t=1}^d(y_{j_v}-y_{b_t})}
      {\prod_{u\ne v}(y_{j_v}-y_{j_u})}.
 \end{split}
\end{equation}
Apply Lemma~\ref{lem:jp-finite} to the ordered factors
\[
 f_\ell(z)=z\quad(1\le\ell\le a) \text{ and }f_{a+t}(z)=z-y_{b_t}\quad(1\le t\le d).
\]
 This gives \eqref{eq:jp-closed}.
 
If $k=0$, every node in \eqref{eq:jp-closed-sum1} lies in $A$, so $E_{r,w}=0$.
For non admissible $w\ne\mathrm{id}$, the vanishing follows from
Lemmas~\ref{lem:jp-basis} and Lemma~\ref{lem:colored}.

To prove positivity, fix a subset $M$ in \eqref{eq:jp-closed}. Factoring
$y_j-y_i=y_i(y_j/y_i-1)$ for $j\le r<i$, we obtain a product
of factors $e^{-\beta}-1$, with $\beta$ a positive root, times
the monomial
\begin{equation*}
 \mu_M=
 \frac{Y_r}{\prod_{i\in L}y_i}
 \prod_{\substack{1\le\ell\le a\\ \ell\notin M}}y_{j_{M_\ell}}
 \prod_{\substack{1\le t\le d\\ a+t\notin M,\ b_t>r}}y_{b_t}.
\end{equation*}
Let $\nu_s$ denote the exponent of $x_s$ in $\mu_M$.
For $s\le r$, we have $\nu_s\ge s-|L\cap[1,s]|\ge0$.
For $s>r$, all $a$ factors from the first block and all tail factors
with $r<b_t\le s$ contribute to $\nu_s$, except at positions in $M$.
Thus
\begin{equation}
\begin{aligned}
 \nu_s
 &\ge r+a+|Z\cap(r,s]|-|M|-|L\cap[1,s]|\\
 &= a+1+|Z\cap[1,s]|-|L\cap[1,s]|\\
 &\ge a+1-|L\setminus R|
  = n+1-|S|\ge0.
\end{aligned}
\end{equation}
Here $|M|=r-|A|-1$, $A=Z\cap[1,r]$, and $0\notin L$.
The last equality follows from $a=n-c-|R|$ and
$|L\cup R|=|S|-c$. Hence $\mu_M\in\mathcal P$, and each
summand of \eqref{eq:jp-closed} belongs to $\mathcal P$.

For $w=\mathrm{id}$, substitute
$(q^{n+1}-1)/(q(1-q))=-\sum_{b=-1}^{n-1}q^b$
into the identity coefficient of \eqref{eq:closedtranslation}. By \eqref{eq:interpnegative}, the $b=-1$
term cancels the constant $1$. By \eqref{eq:interppowers}, the terms with
$0\le b<r-1$ vanish, giving
\[
 (-1)^r E_{r,\mathrm{id}}
 =Y_r\sum_{b=r-1}^{n-1}\mathsf h_{b-r+1}(y_1,\ldots,y_r)
 =Y_r\sum_{s=0}^{n-r}\mathsf h_s(y_1,\ldots,y_r)\in\mathcal P.
\]
This proves \eqref{eq:jp-closed-id}.
\end{proof}

\begin{corollary}\label{cor:largestclosed}
For every nonidentity admissible $w$,
\[
 [D_w]\cg_n=(-1)^{\ell(w)-n}
       \prod_{i\in L}(y_i-1)\prod_{i\in[1,n]\setminus L}y_i,
 \qquad [D_{\id}]\cg_n=(-1)^n Y_n.
\]
\end{corollary}
\begin{proof}
For proper support, $|Z|\le |S|-c\le n-1$.
For full support, choose $j\in L\setminus R$,
so $j\notin Z$. Thus $k=n-|Z|\ge1$.

At $r=n$, we have $A=Z$ and $
 a+d=n-1-|Z|=k-1.$
Hence the unique $M$ consists of every position.
All products over the indices not in $M$ in \eqref{eq:jp-closed} are empty,
leaving
\[
 (-1)^{\ell(w)-n}E_{n,w}
 =
 \frac{Y_n}{\prod_{i\in L}y_i}
 \prod_{i\in L}(y_i-1)
 =
 \prod_{i\in L}(y_i-1)
 \prod_{i\in[1,n]\setminus L}y_i.
\]
Equation \eqref{eq:jp-closed-id} gives the identity coefficient.
\end{proof}

\begin{corollary}[Top-length coefficients]\label{cor:closedtop}
For $j\le r$,
\[
 [T_{\tau_j}]\cg_r=[D_{\tau_j}]\cg_r
 =\prod_{i=1}^n(1-y_i)\prod_{i=r+1}^n(1-y_j/y_i).
\]
For $j>r$ both coefficients vanish, and there are no other coefficients
of length $2n$.
\end{corollary}

\begin{proof}
The data of $\tau_j$ are
$L=[1,n]$ and $R=I\setminus\{j\}$. Hence
\[
 p_{\tau_j}(t)
 =
 \prod_{i=1}^{n}(1-y_i)
 \prod_{\substack{1\le i\le n\\i\ne j}}(1-t/y_i).
\]

Thus $p_{\tau_j}(y_v)=0$ for every interpolation node $v\ne j$. If $j>r$, the sum in \eqref{eq:closedfinite} vanishes. If $j\le r$,
only the term $v=j$ remains. At maximum length the $T$ and $D$ coefficients agree by
\eqref{eq:E}.
\end{proof}

\section{Ordinary and closed Pieri classes in the non-equivariant \texorpdfstring{$K$}{K}-homology}\label{sec:forget}
In this section, we derive the formulas of the coefficients $\phi_0(C_{r,w})$ and $\phi_0(E_{r,w})$. 
For a proper subset $S\subsetneq I$, let $d_S$ be its cyclically decreasing affine Weyl group
element. More precisely, \(d_S\in W_{\mathrm{af}}\) is the element represented by a word containing each \(s_i\), \(i\in S\), exactly once, with \(s_i\) occurring to the left of \(s_{i-1}\) whenever \(i,i-1\in S\). Indices are interpreted modulo \(n+1\). Set \(d_\varnothing=\mathrm{id}\).

For a component $
C=[a,a+1,\ldots,b],$
following the cyclic order, put
\[
d_C=s_b s_{b-1}\cdots s_a,\qquad d_S=\prod_{C:\text{ components of }S}d_C.
\]
Cyclically decreasing elements appear in Lam's construction of affine
Stanley symmetric functions \cite{Lam2006}. Morse's affine set-valued strips give a tableau formulation of the non-equivariant ideal-sheaf Pieri rule \cite{Morse2012}.

For $S\subsetneq I$, define
\[
 I(S):=\{i\in\{1,\ldots,n\}:i-1,i\in S\}.
\]
\begin{theorem}\label{thm:forget}
For $1\le r\le n$, the only nonzero coefficients
$C_{r,w}$ with $\ell(w)=r$ are
\begin{equation}\label{eq:minimal}
 C_{r,d_S}
 =
 \frac{Y_{r-1}}{\prod_{i\in I(S)}y_i},
 \qquad
 S\subsetneq I,\quad |S|=r.
\end{equation}
As a consequence, for $1\le r\le n$ we have
\begin{equation}\label{eq:forget}
 \phi_0(\gh_r)=\sum_{\substack{S\subsetneq I\\|S|=r}}T_{d_S}.
\end{equation}
\end{theorem}
\begin{proof}
Let $m=\ell(w)$. There are three cases. If $m<r$, Corollary~\ref{cor:jp-length} gives $C_{r,w}=0$. If $m>r$, then $C_{r,w}\in \mathfrak m^{m-r}$ and  therefore $\phi_0(C_{r,w})=0$.

Suppose $m=r$. The inequality $m-r\ge|L|$ in the proof of
Corollary~\ref{cor:jp-length} forces $L=\varnothing$.
Full-support data require an $L$-labeled head, so the support is proper.
On each component every internal edge has label $R$ alone.
Its word \eqref{eq:path-word} is therefore decreasing, and the full
word is $d_S$ for a subset of size $r$.
Conversely any such $d_S$ has precisely these data.

Since $Z=I(S)$ and $d=k-1$, only $M=[1,d]$ contributes to
\eqref{eq:jp-ordinary}. The products over $[1,d]\setminus M$ and $L$
are empty, the denominator reduces to $\prod_{i\in I(S)}y_i$, and
the sign is positive since $m=r$. This proves \eqref{eq:minimal}. Under augmentation, each $y_i$ maps
to $1$, giving \eqref{eq:forget}.
\end{proof}

\begin{corollary}\label{cor:closedforget}
For every $1\le r\le n$, we have
\begin{equation}\label{eq:closedforget}
  \phi_0(\cg_r)=\sum_{\substack{S\subsetneq I\\|S|\le r}}T_{d_S}=\sum_{\substack{S\subsetneq I\\|S|\le r}}
(-1)^{r-|S|}
\binom{n-|S|}{r-|S|}
D_{d_S},
\end{equation}
where $d_\varnothing=\id$.
\end{corollary}
\begin{proof}
The first equality follows from taking the sum in \eqref{eq:forget}. Applying $\phi_0$ coefficient wise to the triangular
basis change \eqref{eq:E} gives
\begin{equation}\label{eq:Eaug}
 \phi_0(E_{r,w})=\sum_{a=0}^r
       \sum_{\substack{|S|=a\\w\le d_S}}(-1)^{a-\ell(w)},
\end{equation} 
Every generator occurs exactly once in a reduced expression for \(d_S\). Consequently, the subword criterion gives
\[
w\le d_S
\quad\Longleftrightarrow\quad
w=d_B\ \text{for some }B\subseteq S.
\]
Suppose $w=d_B$ and put $b=|B|$.
Since $I$ contains $n+1$ nodes, there are
 $\binom{n+1-b}{a-b}$
subsets $S$ of size $a$ containing $B$.
Each such subset is proper because $a\le r\le n$.
Therefore
\[
 \phi_0(E_{r,d_B})
 =
 \sum_{a=b}^{r}
 (-1)^{a-b}\binom{n+1-b}{a-b}
 =
 (-1)^{r-b}\binom{n-b}{r-b}.
\]
\end{proof}

\section{Pieri structure constants in two
Schubert bases}\label{sec:pieri}

In this section, we compute the two families $d_{r,v}^{\,u}$ and $\widetilde d_{r,v}^{\,u}$ defined in \eqref{eq_st_con}, prove signed
positivity for $d_{r,v}^{\,u}$, and derive binomial formulas for both non-equivariant specializations.

\begin{proposition}\label{st_con}
   For $u,v\in \Wgr$ and $1\le r\le n$, we have 
\[d_{r,v}^{\,u}=\sum_{\substack{w\in \W\\w\star v=u}}(-1)^{\ell(v)+\ell(w)-\ell(u)}C_{r,w}, \qquad \widetilde d_{r,v}^{\,u}=\sum_{\substack{w\in \W\\(w\star v)^0=u}}E_{r,w}.\]
\end{proposition}
\begin{proof}
We consider the expansion 
$$\gh_r=\sum_{w\in \W}C_{r,w}\,T_w, \qquad \ck_r=\sum_{w\in \W}E_{r,w}\,D_w.$$
We multiply the first equality by $T_v$ from right and the second equality by $D_v$ from right. In the right hand side we use the equality \eqref{eq:demazure}. 
Next, we apply $\kappa$ to $\gh_rT_v$ and $\ck_r D_v$. Using \eqref{eq:moduleprojection} we have $$\kappa(\gh_rT_v)=\gh_r\kappa(T_v)=\gh_r\gh_v, \qquad \kappa(\ck_rD_v)=\ck_r\kappa(D_v)=\ck_r\ck_v$$
Now the proof follows using \eqref{eq:kappa}.
\end{proof}

\begin{theorem}[Signed positivity of ideal-sheaf Pieri structure constants]\label{thm:structure}
For $1\le r\le n$, and $v, u\in\Wg$
\begin{equation}\label{eq:structure}
 (-1)^{\ell(u)-\ell(v)-r}d_{r,v}^{\,u}
 =\sum_{\substack{w\in\W\\w\star v=u}}(-1)^{\ell(w)-r}C_{r,w} \in\mathcal P.
\end{equation}
\end{theorem}

\begin{proof}
Multiply the first identity of
Proposition \ref{st_con} by $(-1)^{\ell(u)-\ell(v)-r}$. Theorem \ref{thm:jp-ordinary} places every resulting summand in $\mathcal P$.
\end{proof}

Proposition~\ref{st_con} expresses $\widetilde{d}_{r,v}^{\,u}$ as the sum of
$E_{r,w}$ over all $w$ satisfying $(w \star v)^0 = u$.
Since the parity of $\ell(w)-r$ can vary within this fiber,
Theorem~\ref{thm:jp-closed} leaves the sign of the sum undetermined. We propose the following conjecture.

\begin{conjecture}\label{conj:closed}
For $1\le r\le n$ and $u,v\in\Wg$,
\[
 (-1)^{\ell(u)-\ell(v)-r}\widetilde d_{r,v}^{\,u}\in\mathcal P.
\]
\end{conjecture}

\begin{remark}
Conjecture~\ref{conj:closed} is the Pieri case of Kumar's conjecture for the structure-sheaf basis~\cite[Conjecture~1.3]{Kumar}, written in the opposite sign convention. Theorem~\ref{thm:structure} establishes the analogous positivity for boundary ideal sheaves.
\end{remark}

\begin{remark}
 Consider finite Weyl group elements
\begin{equation}\label{eq:app-pr}
 w_r:=s_{r-1}\cdots s_1s_\theta\in W,\qquad 1\le r\le n.
\end{equation}
Since $s_0=s_\theta t_{-\theta^\vee}$, we have
$\rho_r=w_rt_{-\theta^\vee}$. Then using \cite[Proposition 2.4]{IKN} and the $K$-theoretic Peterson algebra all quantum product in $QK_T(SL_{n+1}/B)$ (see \cite{MNS1}) with respect to Schubert class corresponding to $w_r$ can be computed explicitly.
\end{remark}

\subsection{Non-equivariant Pieri coefficients in both the bases}
We denote $\bar k_v=\phi_0(k_v)$ and $\bar\ell_v=\phi_0(\ell_v)$. Write
\[
 h_r=\phi_0(k_r),\qquad \widetilde h_r=\phi_0(\ell_r)
       =h_0+\cdots+h_r,\qquad h_0=1.
\]
Consider the products in the non-equivariant $K$-homology
\begin{equation}\label{eq:ns-products}
 h_r\bar k_v=\sum_u\phi_0(d_{r,v}^{\,u})\bar k_u,
 \qquad
 \widetilde h_r\bar\ell_v
       =\sum_u\phi_0(\widetilde d_{r,v}^{\,u})\bar\ell_u.
\end{equation}
\begin{remark}
    Combining Theorem \ref{thm:structure} with Theorem~\ref{thm:forget} gives directly
\begin{equation}\label{eq:ordinarypieri}
 \phi_0(d_{r,v}^{\,u})
 =(-1)^{r+\ell(v)-\ell(u)}
 \#\{S\subset I:|S|=r,\ d_S\star v=u\},
\end{equation}
see also \cite[Corollary~7.6]{LSS}.
\end{remark}

We use the cyclically decreasing elements $d_C$ defined in Section~\ref{sec:forget}.
For $v\in W_{\mathrm{af}}^0$, let
\[
 \mathcal F(v)=\{C\subsetneq I:d_Cv\in W_{\mathrm{af}}^0,
                \ \ell(d_Cv)=\ell(v)+|C|\}.
\]
 
In particular, $\varnothing\in\mathcal F(v)$. For $C\in\mathcal F(v)$,
define
\begin{align}
 a_v(C)&=\#\{i\in I\setminus C:C\cup\{i\}\in\mathcal F(v)\},
       \label{eq:ns-a}\\
 Z_v^-(C)&=\{i\in I\setminus C:C\cup\{i\}\subsetneq I,
                         \ d_{C\cup\{i\}}\star v=d_Cv\},
       \label{eq:ns-zminus}\\
 b_v(C)&=|Z_v^-(C)|.\notag
\end{align}

We use the convention $\binom pq=0$ if $q<0$ or, $q>p$. 
\begin{theorem}
\label{thm:ns-ideal}
Let $u,v\in W_{\mathrm{af}}^0$ and $1\le r\le n$. The coefficients 
$\phi_0(d_{r,v}^{\,u})$ and $\phi_0(\widetilde d_{r,v}^{\,u})$ vanishes unless
$u=d_Cv$ for a unique $C\in\mathcal F(v)$ with $|C|\le r$.
For such $C$,
\begin{align}
 \phi_0(d_{r,v}^{\,d_Cv})&=(-1)^{r-|C|}\binom{b_v(C)}{r-|C|}\label{eq:ns-ideal}\\
 \phi_0(\widetilde d_{r,v}^{\,d_Cv})&=
 \begin{cases}
  1,& \text{if }|C|=r,\\[2pt]
  (-1)^{r-|C|}\displaystyle\binom{a_v(C)-1}{r-|C|}&\text{if }|C|<r.
 \end{cases}\label{eq:ns-closed}
\end{align}

In the second case $a_v(C)\ge1$. Consequently 
$$(-1)^{\ell(u)-\ell(v)-r}\phi_0(d_{r,v}^{\,u})
\in\mathbb Z_{\ge0}, ~(-1)^{\ell(u)-\ell(v)-r}\phi_0(\widetilde d_{r,v}^{\,u})
       \in\mathbb Z_{\ge0}.$$
\end{theorem}

  Theorem \ref{thm:ns-ideal} gives a binomial formula for each coefficient in Takigiku's inclusion--exclusion formula~\cite[Theorem~1.4]{Tak}. The proof of Theorem \ref{thm:ns-ideal} is given in Appendix~\ref{app:nonequivariant}.

\section{Worked examples}\label{sec_ex}

\begin{example}[Type $\mathrm{A}_1$]
  Set $x=e^{-\alpha_1}$. The Pieri classes are
\[
 \gh_1=T_0+T_1-(x-1)T_{01}, \qquad \ck_1=1+\gh_1.
\]
Proposition~\ref{st_con} gives
\[
\gh_1^2=-\kh_{s_0}+\kh_{s_1s_0}-(x-1)\kh_{s_0s_1s_0}.
\]
The signs
are $(-1)^{\ell(u)-2}$, and the signed weights are respectively
$1,1,x-1$. At $x=1$ we have
$\phi_0(\gh_1)^2=-\phi_0(\kh_{s_0})+\phi_0(\kh_{s_1s_0})$. Consequently
\[
 \cg_1^2=x\ck_{s_1s_0}+(1-x)\ck_{s_0s_1s_0}.
\]
The two signed coefficients are $x$ and $x-1$. At $x=1$ this becomes
$\phi_0(\cg_1)^2=\phi_0(\ck_{s_1s_0})$.
\end{example}

\begin{example}[Type $\mathrm{A}_2$]\label{eg_A2}
We denote $x=x_1=e^{-\alpha_1}$, $y=x_2=e^{-\alpha_2}$, $a=y-1$, $b=x-1$, $c=xy-1$. The Pieri classes are
\begin{align*}
 \gh_1={}&T_0+T_1+T_2-c(T_{01}+T_{02})-aT_{12}-bT_{21}+acT_{012}+bcT_{021}+abT_{121}-abcT_{0121},\\[3pt]
 \gh_2={}&xyT_{02}+T_{10}+xT_{21}-cT_{010}-cxT_{021}-aT_{102}-axT_{121}+acxT_{0121}+acT_{1012}.
\end{align*}

Using $\cg_2=1+\gh_1+\gh_2$, we obtain
\begin{align*}
 \cg_2={}&1+T_0+T_1+T_2-cT_{01}+T_{02}+T_{10}-aT_{12}+T_{21}\\
 &-c(T_{010}+T_{021})+acT_{012}-a(T_{102}+T_{121})+ac(T_{0102}+T_{0121}).
\end{align*}
In the $\{D_w\}$ basis the same class is
\begin{align*}
 \cg_2={}&xy^2(D_{\id}-D_0-D_1-D_2)+ycD_{01}
          +xy^2(D_{02}+D_{10}+D_{21})+xyaD_{12}\\
 &-yc(D_{010}+D_{021})-acD_{012}
          -xya(D_{102}+D_{121})+ac(D_{0102}+D_{0121}).
\end{align*}
\end{example}

\begin{example}
In $\mathrm{A}_3$, choose $L=\{1,3\}$ and $R=\{0,2,3\}$.
Head $3$ has label $B$ and head $1$ is the pivot. Corresponding admissible element is
$w=s_0s_2s_3s_2s_1$, of length $5$, and
\[
 p_w(t)=(1-y_1)(1-y_3)(1-t/y_2)(1-t/y_3).
\]
 Theorem \ref{thm:jp-ordinary} gives
\begin{align*}
 C_{1,w}&=(y_1-1)(y_3-1)(y_1/y_2-1)(y_1/y_3-1),\\
 C_{2,w}&=-\frac{y_1}{y_2}(y_1-1)(y_3-1)(y_1/y_3-1),\\
 C_{3,w}&=\frac{y_1}{y_3}(y_1-1)(y_3-1).
\end{align*}
Each signed
coefficient is a monomial times a product of factors $e^{-\beta}-1$ witth $\beta>0$, and all three
specialize to zero because $\ell(w)>r$.

 In $\mathrm{A}_3$, consider $w=s_0s_3s_2$ and $r=2$,. Its data are
$L=\varnothing$, $R=\{0,3\}$, $c=1$, so
$a=0$, $h=1$, $A=\varnothing$, $J=(1,2)$.
The index list in \eqref{eq:jp-list} is $(b_1,b_2)=(0,3)$, corresponding to the characters $(1,y_3)$.
Since $\ell(w)-r=1$, Theorem \ref{thm:jp-closed} gives
\[
 -E_{2,032}=y_1y_2\bigl((y_2-y_3)+(y_1-1)\bigr)
 =y_1y_2\bigl(y_3(y_2/y_3-1)+(y_1-1)\bigr)\in\mathcal P.
\]
\end{example}

\begin{example}
In $\mathrm{A}_5$, take $S=\{0,2,4\}$, whose three components are singletons. Then $w=s_0s_2s_4$, $m=3$, $L=R=\varnothing$, $c=3$, $h=2$.
For $r=2$, $A=\varnothing$, $J=(1,2)$, and the index list is
$(b_1,b_2)=(0,0)$, with $y_{b_1}=y_{b_2}=1$. Here $d=2$
although $\ell(w)=3$. The common
prefactor is $y_1$, and the two complete terms for $M=\{1\}$
and $M=\{2\}$ are respectively $y_1(y_2-1)$ and $y_1(y_1-1)$. Thus
\[
 C_{2,024}=-y_1\bigl((y_1-1)+(y_2-1)\bigr).
\]
For $r= 1$ and 3, \eqref{eq:jp-ordinary} gives $C_{1,024}=(y_1-1)^2$ and $C_{3,024}=y_1y_2$, respectively.
\end{example}

\appendix

\section{Proof of Theorem \ref{thm:ns-ideal}}\label{app:nonequivariant}

Fix $v\in W_{\mathrm{af}}^0$ and $B\subsetneq I$, define
\[
  \Psi_v(B)=d_B\star v,\qquad \Phi_v(B)=(d_B\star v)^0.
\]
For $C\in\mathcal F(v)$, let
\[
  Z_v^+(C)=\{i\in I\setminus C \mid C\cup\{i\}\subsetneq I,\ 
                \Phi_v(C\cup\{i\})=d_Cv\}.
\]
Recall that $Z_v^-(C)$ is defined in the same way using $\Psi_v$.

\begin{lemma}\label{lem:non-equivariant-pieri-fibers}
Every Grassmannian value of $\Psi_v$ and every value of $\Phi_v$ is
$d_Cv$ for a unique $C\in\mathcal F(v)$. For $B\subsetneq I$,
\[
\begin{aligned}
  \Psi_v(B)=d_Cv &\iff C\subseteq B\subseteq C\cup Z_v^-(C),\\
  \Phi_v(B)=d_Cv &\iff C\subseteq B\subseteq C\cup Z_v^+(C).
\end{aligned}
\]
If $|C|<n$, then $a_v(C)\ge1$ and
$|Z_v^+(C)|=n+1-|C|-a_v(C)$.
\end{lemma}

\begin{proof}
Let $w_0$ be the longest element of $W$. Since $z\star w_0=z^0w_0$,
we have $$d_B\star(vw_0)=\Phi_v(B)w_0.$$ Thus it suffices to consider
$d_B\star(vq)$ for $q\in\{e,w_0\}$. Read a cyclically decreasing word
for $d_B$ from right to left, and let $C$ be the set of indices at which
the length increases. The corresponding subword represents $d_C$, so
if the endpoint is $uq$ with $u\in W_{\mathrm{af}}^0$, then
$uq=d_Cvq$ and $\ell(uq)-\ell(vq)=|C|$. By parabolic length additivity,
$u=d_Cv$ and $\ell(u)-\ell(v)=|C|$. Hence $C\in\mathcal F(v)$, and
$C=\operatorname{Supp}(uv^{-1})$ proves uniqueness. 

Now fix
$C\in\mathcal F(v)$ with $C\subseteq B$. If $d_B\star(vq)=d_Cvq$,
every letter indexed by $B\setminus C$ fixes the current element, so
retaining any one of these letters together with the $C$-letters
preserves the endpoint. Conversely, if each such restriction has
endpoint $d_Cvq$, uniqueness shows that each extra letter fixes the
element obtained by applying the $C$-letters to its right to $vq$.
Induction from right to left then gives $d_B\star(vq)=d_Cvq$.
Taking $q=e$ and $q=w_0$ gives the two fiber descriptions.

Finally, suppose $|C|<n$. The weak-strip family
$\mathcal F(v)=Z'_{vw_0,+}$ has a greatest member of size $n$ and
admits saturated chains under inclusion
\cite[Proposition~4.12(2), Remark~4.13, and Corollary~4.17]{Tak}.
Thus $C$ has a one-element extension in $\mathcal F(v)$, so
$a_v(C)\ge1$. For each $i\in I\setminus C$, Bruhat monotonicity gives
$d_Cv\le\Phi_v(C\cup\{i\})$. Since $\ell(d_Cv)=\ell(v)+|C|$,
the inequality is strict precisely when all $|C|+1$ steps increase
length, equivalently when $C\cup\{i\}\in\mathcal F(v)$. Equality
holds precisely when $i\in Z_v^+(C)$. Counting the $n+1-|C|$
possible indices yields $|Z_v^+(C)|=n+1-|C|-a_v(C)$.
\end{proof}

\begin{proof}[Proof of Theorem \ref{thm:ns-ideal}]
By Proposition~\ref{st_con} and Corollary~\ref{cor:closedforget}, only the endpoints described in
Lemma~\ref{lem:non-equivariant-pieri-fibers} with $|C|\le r$ can
contribute. Fix $u=d_Cv$. In the ideal-sheaf basis, \eqref{eq:ordinarypieri} counts the
sets $B=C\cup J$ with $J\subseteq Z_v^-(C)$ and $|J|=r-|C|$, each
with sign $(-1)^{r-|C|}$. Since $|Z_v^-(C)|=b_v(C)$, this gives \eqref{eq:ns-ideal}.

In the structure-sheaf basis, a set $B$ of size $s\le r$ contributes
$(-1)^{r-s}\binom{n-s}{r-s}$ by \eqref{eq:closedforget}. If $|C|=r$, only $B=C$
contributes, with coefficient $1$. Suppose $|C|<r$ and put
$M=|Z_v^+(C)|=n+1-|C|-a_v(C)$. Since $a_v(C)\ge1$, we have
$M\le n-|C|$. Grouping the sets $B=C\cup J$ by $j=|J|$ gives
\[
\begin{aligned}
 \phi_0(\widetilde d^{\,u}_{r,v})
 &= [z^{r-|C|}]\sum_{j=0}^{M}\binom{M}{j}z^j(1-z)^{n-|C|-j}\\
 &= [z^{r-|C|}](1-z)^{n-|C|-M}
  =(-1)^{r-|C|}\binom{a_v(C)-1}{r-|C|}.
\end{aligned}
\]
 This proves \eqref{eq:ns-closed}, and
both sign assertions follow from $\ell(u)-\ell(v)=|C|$.
\end{proof}

 \noindent  {\bf Acknowledgment.} The author is especially grateful to Takeshi Ikeda for many insightful discussions throughout this project. This work was supported by JSPS KAKENHI Grant Numbers 24KF0258, 25KF0074. The author was supported by a JSPS Postdoctoral Fellowship for Research in Japan, during which this research was carried out. The author thanks Waseda University for providing a supportive research environment.


\begin{thebibliography}{99}

\bibitem{AJS}
H.~H. Andersen, J.~C. Jantzen, and W.~Soergel,
Representations of quantum groups at a $p$-th root of unity and of semisimple
groups in characteristic $p$: independence of $p$,
\emph{Ast\'erisque}, no.~220 (1994). MR1272539.


\bibitem{AGM2011}
D.~Anderson, S.~Griffeth, and E.~Miller,
Positivity and Kleiman transversality in equivariant $K$-theory of homogeneous spaces,
\emph{J. Eur. Math. Soc. (JEMS)} \textbf{13} (2011), no.~1, 57--84. MR2735076.

\bibitem{BK2017}
S.~Baldwin and S.~Kumar,
Positivity in $T$-equivariant $K$-theory of flag varieties associated to
Kac--Moody groups II,
\emph{Represent. Theory} \textbf{21} (2017), 35--60. MR3627147.

\bibitem{Bil}
S.~C. Billey,
Kostant polynomials and the cohomology ring for $G/B$,
\emph{Duke Math. J.} \textbf{96} (1999), no.~1, 205--224. MR1663931.

\bibitem{BB}
A.~Bj\"orner and F.~Brenti,
\emph{Combinatorics of Coxeter Groups},
Graduate Texts in Mathematics, vol.~231, Springer, Berlin, 2005. MR2133266.

\bibitem{CL}
C.~H. Chow and N.~C. Leung,
Quantum $K$-theory of $G/P$ and $K$-homology of affine Grassmannian,
arXiv:2201.12951 (2022).

\bibitem{deBoor}
C. de Boor, Divided differences,
\emph{Surv. Approx. Theory} \textbf{1} (2005), 46--69.

\bibitem{FK1}
S.~Fomin and A.~N. Kirillov,
Yang--Baxter equation, symmetric functions and Grothendieck polynomials,
arXiv:hep-th/9306005 (1993).

\bibitem{FK2}
S.~Fomin and A.~N. Kirillov,
Grothendieck polynomials and the Yang--Baxter equation,
in \emph{Formal Power Series and Algebraic Combinatorics /
S\'eries formelles et combinatoire alg\'ebrique} (Piscataway, NJ, 1994),
DIMACS, Piscataway, NJ, 1994, pp.~183--189. MR2307216.

\bibitem{IIM}
T.~Ikeda, S.~Iwao, and T.~Maeno,
Peterson isomorphism in $K$-theory and relativistic Toda lattice,
\emph{Int. Math. Res. Not. IMRN} \textbf{2020}, no.~19, 6421--6462. MR4165480.

\bibitem{IKN}
T.~Ikeda, T.~Kouno, and S.~Naito,
Seidel product formula in equivariant quantum $K$-theory of flag varieties,
arXiv:2602.24050v1 (2026).

\bibitem{ISY}
T.~Ikeda, M.~Shimozono, and K.~Yamaguchi,
Equivariant $K$-homology of affine Grassmannian and $K$-theoretic double
$k$-Schur functions,
\emph{Adv. Math.} \textbf{492} (2026), article~110894, 57~pp. MR5039968.

\bibitem{Kato}
S.~Kato,
Loop structure on equivariant $K$-theory of semi-infinite flag manifolds,
\emph{Ann. of Math.} (2) \textbf{202} (2025), no.~3, 1001--1075. MR4979913.

\bibitem{KK}
B.~Kostant and S.~Kumar,
The nil Hecke ring and cohomology of $G/P$ for a Kac--Moody group $G$,
\emph{Adv. in Math.} \textbf{62} (1986), no.~3, 187--237. MR0866159.

\bibitem{KK1}
B.~Kostant and S.~Kumar,
$T$-equivariant $K$-theory of generalized flag varieties,
\emph{J. Differential Geom.} \textbf{32} (1990), no.~2, 549--603. MR1072919.

\bibitem{Kumar}
S.~Kumar,
Conjectural positivity for Pontryagin product in equivariant $K$-theory of loop groups,
\emph{J. Pure Appl. Algebra} \textbf{230} (2026), no.~7, article~108274, 37~pp. MR5070875.

\bibitem{Lam2006}
T.~Lam,
Affine Stanley symmetric functions,
\emph{Amer. J. Math.} \textbf{128} (2006), no.~6, 1553--1586. MR2275911.

\bibitem{LLMSS}
T.~Lam, L.~Lapointe, J.~Morse, A.~Schilling, M.~Shimozono, and M.~Zabrocki,
\emph{$k$-Schur Functions and Affine Schubert Calculus},
Fields Institute Monographs, vol.~33, Springer, New York, 2014. MR3379711.

\bibitem{LLMS}
T.~Lam, L.~Lapointe, J.~Morse, and M.~Shimozono,
Affine insertion and Pieri rules for the affine Grassmannian,
\emph{Mem. Amer. Math. Soc.} \textbf{208} (2010), no.~977, xii+82~pp.

\bibitem{LLMSC}
T.~Lam, C.~Li, L.~C. Mihalcea, and M.~Shimozono,
A conjectural Peterson isomorphism in $K$-theory,
\emph{J. Algebra} \textbf{513} (2018), 326--343. MR3849889.

\bibitem{LSS}
T.~Lam, A.~Schilling, and M.~Shimozono,
$K$-theory Schubert calculus of the affine Grassmannian,
\emph{Compos. Math.} \textbf{146} (2010), no.~4, 811--852. MR2660675.

\bibitem{LSSC}
T.~Lam, A.~Schilling, and M.~Shimozono,
Schubert polynomials for the affine Grassmannian of the symplectic group,
\emph{Math. Z.} \textbf{264} (2010), no.~4, 765--811. MR2593294.

\bibitem{LS-Q}
T.~Lam and M.~Shimozono,
Quantum cohomology of $G/P$ and homology of affine Grassmannian,
\emph{Acta Math.} \textbf{204} (2010), no.~1, 49--90. MR2600433.

\bibitem{LS}
T.~Lam and M.~Shimozono,
Equivariant Pieri rule for the homology of the affine Grassmannian,
\emph{J. Algebraic Combin.} \textbf{36} (2012), no.~4, 623--648. MR2984159.

\bibitem{LLT}
A.~Lascoux, B.~Leclerc, and J.-Y.~Thibon,
Flag varieties and the Yang--Baxter equation,
\emph{Lett. Math. Phys.} \textbf{40} (1997), no.~1, 75--90. MR1445968.

\bibitem{MNS1}
T.~Maeno, S.~Naito, and D.~Sagaki,
A presentation of the torus-equivariant quantum $K$-theory ring of flag
manifolds of type~$A$, Part~I: The defining ideal,
\emph{J. Lond. Math. Soc.} (2) \textbf{111} (2025), no.~3, article~e70095.

\bibitem{Morse2012}
J.~Morse,
Combinatorics of the $K$-theory of affine Grassmannians,
\emph{Adv. Math.} \textbf{229} (2012), no.~5, 2950--2984. MR2889152.

\bibitem{Pet}
D.~Peterson,
\emph{Quantum cohomology of $G/P$},
lecture notes, Massachusetts Institute of Technology, Cambridge, MA, spring 1997.

\bibitem{Pon}
S.~Pon,
Affine Stanley symmetric functions for classical types,
\emph{J. Algebraic Combin.} \textbf{36} (2012), no.~4, 595--622. MR2984158.

\bibitem{Tak}
M.~Takigiku,
A Pieri formula and a factorization formula for sums of $K$-theoretic $k$-Schur functions,
\emph{Algebr. Comb.} \textbf{2} (2019), no.~4, 447--480. MR3997506.


\end{thebibliography}
\end{document}